\documentclass[11pt]{article}
\usepackage{latexsym,amsmath,amscd,amssymb, graphics}
\usepackage{amsmath,amssymb,mathrsfs,framed,esint,slashed,color}
\usepackage{amsthm}
\usepackage{enumerate}
\usepackage{hyperref}
\hypersetup{
colorlinks = true, %Colours links instead of ugly boxes
urlcolor = blue, %Colour for external hyperlinks
linkcolor = blue, %Colour of internal links
citecolor = red %Colour of citations
}
\usepackage{enumitem}
\usepackage[margin=1in]{geometry}
\usepackage{stmaryrd}
\usepackage{graphicx}
\usepackage[all]{xy}
\usepackage[makeroom]{cancel}
\usepackage{empheq}

\usepackage{framed}

\usepackage{pgfplots, pgfplotstable, booktabs, colortbl, array,multirow}
\usetikzlibrary{calc}
\usepgfplotslibrary{groupplots}
\pgfplotsset{compat=newest}

\theoremstyle{plain}

\theoremstyle{definition}

\newtheorem{remark}{Remark}[section]

\DeclareMathOperator{\dv}{div}
\DeclareMathOperator{\curl}{curl}

\usepackage{verbatim}

\begin{document}

\title{Structure-preserving and thermodynamically consistent finite element discretization for visco-resistive
MHD with thermoelectric effect}

\author{Evan S. Gawlik\thanks{\noindent 
Department of Mathematics and Computer Science,  Santa Clara University, \href{egawlik@scu.edu}{egawlik@scu.edu}},\;
Fran\c{c}ois Gay-Balmaz\thanks{\noindent Division of Mathematical Sciences, Nanyang Technological University, Singapore \href{francois.gb@ntu.edu.sg}{francois.gb@ntu.edu.sg}}
, and Bastien Manach-P\'erennou\thanks{\noindent Division of Mathematical Sciences, Nanyang Technological University, Singapore
\href{bastien.manachp@ntu.edu.sg}{bastien.manachp@ntu.edu.sg}}
}

\date{}

\maketitle

\begin{abstract} 
We present a structure-preserving and thermodynamically consistent numerical scheme for classical magnetohydrodynamics, incorporating viscosity, magnetic resistivity, heat transfer, and thermoelectric effect.  The governing equations are shown to be derived from a generalized Hamilton's principle, with the resulting weak formulation being mimicked at the discrete level. The resulting numerical method conserves mass and energy, satisfies Gauss' magnetic law and magnetic helicity balance, and adheres to the Second Law of Thermodynamics, all at the fully discrete level. It is shown to perform well on magnetic Rayleigh-B\'enard convection.
\end{abstract}

%\tableofcontents

\section{Introduction}

\textbf{Dissipative magnetohydrodynamics.} Plasma flows are at the core of numerous applications, including but not limited to: astrophysics, Magnetic Confinement Fusion \cite{Chen2016,Stacey2005} and Inertial Confinement Fusion \cite{AtMtV2004,Drake2006}. Such flows display complex behaviors because of thermodynamics, interaction with the electromagnetic field, potentially highly compressible effects, and numerous and stiff irreversible processes. When characteristic velocities are negligible compared to the speed of light and the plasma is highly collisional, the flow can be reasonably described by Magnetohydrodynamics (MHD). Even for flows where the plasma is less collisional (i.e. for Magnetic Confinement Fusion), these equations are often used as a first approximation. The system under interest is
\begin{subequations}
\label{eq:full_MHD}
\begin{align}
\rho \left( \partial_t u + u \cdot \nabla u \right) + \nabla p & = \frac{1}{\mu_0} \operatorname{curl} B \times B +\nabla \cdot \sigma, \\
\partial_t B + \operatorname{curl} ( B \times u) & = - \frac{1}{\mu_0} \operatorname{curl} (\nu \operatorname{curl} B) - \operatorname{curl} (\alpha \nabla T), \\
\partial_t \rho + \nabla \cdot (\rho u) & = 0, \\
T \left( \partial_t s + \nabla \cdot (s u) \right) & = \sigma : \nabla u + \nabla \cdot (\kappa \nabla T) + \frac{\nu}{\mu_0^2} \operatorname{curl} B \cdot \operatorname{curl} B - \frac{T}{\mu_0} \nabla \cdot (\alpha \operatorname{curl} B).
\end{align}
\end{subequations}
The notations are standard: $u$ is the velocity, $B$ the magnetic field, $\rho$ the mass density, $s$ the entropy density, $p$ the pressure and $T$ the temperature; $\mu_0$ is the vacuum permeability, $\sigma$ the viscous tensor, $\nu$ the resistivity coefficient and $\kappa$ the thermal conductivity. Finally, $\alpha$ denotes the thermoelectric coefficient. Equations \eqref{eq:full_MHD} will be supplemented with the boundary conditions
\begin{equation}
u = B \cdot n = \operatorname{curl} B \times n = 0,
\end{equation}
with $n$ the outward normal vector to the boundary, together with a thermal boundary condition of either Dirichlet or Neumann type. Extension of these equations to include thermoelectric effect will be also considered.
When dissipation (i.e. viscosity, resistivity and heat exchange) is removed, the so-called ideal MHD is recovered. It is well-known that ideal MHD equations can be derived from Hamilton's principle \cite{Newcomb1962}. Because dissipation is by nature irreversible and leads to an entropy production, it falls in the realm of non-equilibrium thermodynamics. An extension of Hamilton's principle for non-equilibrium thermodynamics was developed in \cite{GayBalmaz2017a, GayBalmaz2017b}. It allows to derive equations automatically satisfying the two first laws of thermodynamics; namely conservation of energy and non-destruction of entropy. The above equations \eqref{eq:full_MHD}, including all irreversible processes, can be obtained from it, as will be shown shortly. This structure helps navigate through the otherwise complicated and intricate phenomena. More importantly, it will here serve as a guideline for deriving a numerical method that preserves key properties of the flow at the discrete level. Dissipative MHD has also been studied through extensions of Poisson brackets, including the metriplectic approach; see \cite{MaTa2012} and \cite{Coquinot2020}.

\paragraph{Structure-preserving discretization.} In the literature, little work tackles the numerical approximation of system \eqref{eq:full_MHD}. Most of the time, incompressible flow is assumed, irreversible processes are neglected or the entropy equation is removed in favor of a barotropic closure. Structure preserving methods for incompressible and ideal MHD with constant density include \cite{LiWa2001}, \cite{GaMuPaMaDe2011,HuMaXu2017,KrMa2017,HiLoMaZh2018,HuXu2019,HuLeXu2021}, while the variable density case was treated in \cite{GaGB2022}. These methods have succeeded in preserving at the discrete levels several invariants and
constraints of the continuous system. For instance, in \cite{GaGB2022},
a finite element method was proposed which preserves energy, cross-helicity (when
the fluid density is constant), magnetic helicity, mass, total squared density, pointwise incompressibility and Gauss' magnetic law to machine precision, both at the spatially and temporally discrete levels. Structure-preserving methods for compressible MHD were developed in \cite{Gawlik2021} and  \cite{Carlier2025}. While the process of resistivity and viscosity are included in \cite{Gawlik2021}, they are not treated as irreversible processes but as dissipative effects with no impact on the entropy equation.

Preserving the structure and physical laws becomes all the more difficult when considering the full system \eqref{eq:full_MHD} and its irreversible processes.
The present numerical method is derived from mimicking the weak form of the equations obtained from the extension of Hamilton's principle to non-equilibrium thermodynamics developed in \cite{GayBalmaz2017b}. The resulting scheme preserves mass, energy, and the pointwise Gauss' magnetic law, while consistently reproducing the magnetic helicity balance and guaranteeing positive entropy production, after both spatial and temporal discretization. Overall, using a discrete variational principle has major benefits on the accuracy, stability and long-term behavior of numerical approximation. In particular, the variational approach ensures that discrete entropy production is effectively driven by physical dissipation and does not suffer from large spurious numerical traces. It allows for accurate long-term behavior and is critical for systems where physical stability depends on a subtle balance between the different sources of dissipation.

\paragraph{Organization of the paper.} In Section \ref{sec_2}, the variational formulation of the equations is introduced. The derivation of the reversible equations from Hamilton's principle is first presented and then extended to include vicosity, resistivity, heat exchange, and thermoelectric effect.
This section provides the foundation for the numerical method developed in Section \ref{sec_3}, where the equations are discretized in space using a discrete version of the variational formulation. The resulting spatially discrete system is then integrated in time using a scheme that maintains all the structure-preserving properties of the spatial discretization. Finally, several test cases are considered in Section \ref{sec_4} to showcase both the stability and accuracy of the scheme.

\section{Variational formulation}\label{sec_2}

\subsection{Ideal MHD}

\paragraph{Lagrangian description.} It is well known that in the absence of irreversible processes, the equations of motion for magnetohydrodynamics follow from Hamilton's principle in the Lagrangian description, from which the Eulerian variational formulation is deduced by using the relabeling symmetry.

Assume that the fluid moves in a compact domain $ \Omega \subset \mathbb{R} ^3$, with a smooth boundary. We denote by $ \operatorname{Diff}( \Omega )$ the group of diffeomorphisms of $\Omega $ and by $ \varphi :[t_0,t_1] \rightarrow \operatorname{Diff}( \Omega )$ the flow of the fluid. We write $ \varrho\, {\rm d} x:[t_0,t_1] \rightarrow \Lambda ^3( \Omega )$, $ S\, {\rm d} x:[t_0,t_1] \rightarrow \Lambda ^3( \Omega )$ and $ \mathcal{B}  \cdot {\rm d} A:[t_0,t_1] \rightarrow  \Lambda ^2( \Omega )$ the mass density, entropy density, and magnetic induction in the Lagrangian description. They are related to the Eulerian mass density $ \rho \, {\rm d} x: [t_0,t_1] \rightarrow \Lambda ^3( \Omega )$, Eulerian entropy $s  {\rm d} x: [t_0,t_1] \rightarrow \Lambda ^3( \Omega )$ and Eulerian magnetic induction $B \cdot {\rm d} a:[t_0,t_1] \rightarrow  \Lambda ^2( \Omega )$ by the push-foward operations of $3$-forms and $2$-forms
\begin{subequations}
\label{eq:push_forward} 
\begin{align}
\rho (t) {\rm d} x & = \varphi (t)_ *( \varrho (t) {\rm d} x), \\ 
s (t) {\rm d} x & = \varphi (t)_ *( S (t) {\rm d} x),  \\
B  (t) \cdot {\rm d} a & = \varphi (t)_ *( \mathcal{B}  (t) \cdot {\rm d} A).
\end{align}
\end{subequations}
In the reversible case, as a consequence of mass and entropy conservation, and of the frozen-in property of magnetic induction, one has $\varrho (t)= \varrho _0$, $S(t)=S_0$ and $ \mathcal{B} (t)= \mathcal{B} _0$ in the Lagrangian description. In this case, the MHD motion follows from Hamilton's principle
\begin{equation}\label{eq:HP_MHD}
\delta \int_{t_0}^{t_1} L( \varphi , \dot  \varphi , \varrho _0, S_0, \mathcal{B} _0) {\rm d} t=0 ,
\end{equation} 
for arbitrary variations $ \delta \varphi $ vanishing at $t=t_0,t_1$.

\paragraph{Eulerian description.} From the relabelling symmetry, the Lagrangian can be written in terms of the Eulerian variables \eqref{eq:push_forward} as
\begin{equation}
L( \varphi , \dot  \varphi , \varrho _0, S_0, \mathcal{B} _0)= \ell(u, \rho  , s, B),
\end{equation}
with $u= \dot  \varphi \circ \varphi ^{-1} $ the Eulerian velocity. The Eulerian description naturally involves covariant time derivatives and covariant variations, which are denoted as: 
\begin{subequations}
\label{eq:covariant_deriv}
\begin{align} 
D_t \gamma & = \partial_t \gamma + u \cdot \nabla \gamma, & 
D_\delta \gamma & = \delta \gamma + v \cdot \nabla \gamma, \\
\mathfrak{D}_t Z & = \partial_t Z + \operatorname{curl}Z \times u + \nabla (Z \cdot u), & 
\mathfrak{D}_\delta Z & = \delta Z + \operatorname{curl}Z \times v + \nabla (Z \cdot v), \\
\overline{\mathfrak{D}}_t B & = \partial_t B+ \operatorname{curl}(B \times u) + u \operatorname{div}B, &
\overline{\mathfrak{D}}_\delta B&= \delta B+ \operatorname{curl}(B \times v) + v \operatorname{div}B, \\
\overline{D}_t s & = \partial _t s + \operatorname{div}(su), & \overline{D}_\delta s & = \delta s + \operatorname{div}(sv),
\end{align} 
\end{subequations}
with $u$ the Eulerian velocity and $v$ the Eulerian variation.
These four expressions are essentially the same mathematical expression, but applied to forms of different degree. In terms of the Lie derivative of $k$-forms, $k \in \{0,1,2,3\}$, they respectively read 
\begin{subequations}
\begin{align}
D_t \gamma & = \partial_t \gamma + \pounds_u \gamma, &
D_\delta \gamma & = \delta \gamma + \pounds_v \gamma, \\
\mathfrak{D}_t Z \cdot {\rm d}x & = \partial_t Z \cdot {\rm d} x + \pounds_u (Z \cdot {\rm d} x), &
\mathfrak{D}_\delta Z \cdot {\rm d}x & = \delta Z \cdot {\rm d} x + \pounds_v (Z \cdot {\rm d} x), \\
\overline{\mathfrak{D}}_t B \cdot {\rm d}a & = \partial_t B \cdot {\rm d} a + \pounds_u (B \cdot {\rm d} a), &
\overline{\mathfrak{D}}_\delta B \cdot {\rm d}a & = \delta B \cdot {\rm d} a + \pounds_v (B \cdot {\rm d} a), \\
\overline{D}_t \rho {\rm d}x & = \partial_t \rho {\rm d} x + \pounds_u (\rho {\rm d} x), &
\overline{D}_\delta \rho {\rm d}x & = \delta \rho {\rm d} x + \pounds_v (\rho {\rm d} x).
\end{align}
\end{subequations}
The $0$-form $\gamma$ and $1$-form $Z$ will be used later in the irreversible case.
The reduced Hamilton's principle \eqref{eq:HP_MHD} then becomes
\begin{equation}
\label{eq:EP} 
\delta \int_{t_0}^{t_1}\ell(u, \rho  , s, B) {\rm d} t=0
\end{equation} 
for variations of the form
\begin{equation}
\label{eq:EP_var}
\delta u = \partial _t v +\pounds _u v, \quad \overline{D}_\delta \rho  = 0, \quad \overline{D}_\delta s = 0,  \quad \overline{\mathfrak{D}}_\delta B = 0,
\end{equation} 
where $v:[t_0,t_1] \rightarrow \mathfrak{X} ( \Omega )$ is an arbitrary time dependent vector field parallel to the boundary, with $v(t_0)=v(t_1)=0$. A direct application of \eqref{eq:EP}--\eqref{eq:EP_var} yields the fluid momentum equations in the form
\begin{equation}\label{eq:EP_weak_form} 
\left\langle \partial _t \frac{\delta \ell}{\delta u} , v \right\rangle + a \left( \frac{\delta \ell}{\delta u}, u, v \right) +  b ^1  \left( \frac{\delta \ell}{\delta \rho  }, \rho  , v \right)+  b^1 \left( \frac{\delta \ell}{\delta s  }, s  , v \right)+ b^2 \left( \frac{\delta \ell}{\delta B}, B, v \right)=0,
\end{equation}
for all $ v $ with $ v \cdot n=0$, with trilinear forms
\begin{subequations}
\begin{align}
a(w,u,v)&= - \int_ \Omega w \cdot [u,v] {\rm d} x \label{eq:awuv}\\
b^1( \sigma , \rho  , v)&=- \int_ \Omega \rho  \nabla \sigma \cdot v {\rm d} x \label{eq:sigmarhov}\\
b^2(C, B, v)&= \int_ \Omega C \cdot \operatorname{curl}(B \times v) {\rm d} x \label{eq:cCBv}.
\end{align}
\end{subequations}
The equations for $ \rho  $, $s$, and $B$ follow from their definition in \eqref{eq:push_forward} with $ \varrho (t)= \varrho _0$, $S(t)=S_0$, and $ \mathcal{B} (t)= \mathcal{B} _0$, which are expressed in terms of $b^1$ and $b^2$ as
\begin{subequations}
\label{eq:transport_weak_form}
\begin{align}
\left\langle \partial _t \rho  , \sigma \right\rangle +b^1( \sigma , \rho  , u) & =0, \quad \forall \,\sigma\\
\left\langle \partial _t s  , \sigma \right\rangle +b^1( \sigma , s  , u) & =0, \quad \forall \,\sigma\\
\left\langle \partial _t B, C \right\rangle +b^2(C,B,u) & =0, \quad \forall \,C, \ C \cdot n=0.
\end{align} 
\end{subequations}
The strong form of \eqref{eq:EP_weak_form}-\eqref{eq:transport_weak_form} is the Euler-Poincar\'e equation supplemented by conservation of mass and entropy, and the ideal induction law
\begin{subequations}
\begin{align}
\partial _t \frac{\delta \ell}{\delta u} + \pounds _u \frac{\delta \ell}{\delta u} & = \rho \nabla \frac{\delta \ell}{\delta \rho  }+s \nabla \frac{\delta \ell}{\delta s} + B \times \operatorname{curl} \frac{\delta \ell}{\delta B}, \\
\label{eq:strong_form_rho}
\overline{D}_t \rho & = 0, \\
\label{eq:strong_form_s}
\overline{D}_t s & = 0, \\
\label{eq:strong_form_B}
\overline{\mathfrak{D}}_t B & = 0.
\end{align}
\end{subequations}

\subsection{Visco-resistive MHD with heat conduction and thermoelectric effect}

\paragraph{Lagrangian description.} For visco-resistive MHD with heat conduction and thermoelectric effect, besides the Lagrangian, one also needs to specify the phenomenological expressions of the viscous stress tensor, entropy flux, and resistive flux, denoted $ P $, $J_S$, and $ J_ \mathcal{B} $ in the Lagrangian description. The nature of these fluxes is as follows: $P$ is a (1,1) two point tensor field density covering $ \varphi $, \cite{Marsden1994}, $J_S \otimes {\rm d} ^3X$ is a vector field density, and $ J_ \mathcal{B}  \cdot {\rm d} X$ is a one-form. As opposed to the reversible case treated earlier, in the visco-resistive case the no-slip boundary condition for the velocity has to be enforced, which implies the use of the subgroup $ \operatorname{Diff}_0( \Omega ) \subset \operatorname{Diff}( \Omega )$ of diffeomorphisms which fix the boundary pointwise. 

We formulate an extension of Hamilton's principle \eqref{eq:HP_MHD} to visco-resistive MHD with heat transfer by extending the variational approach of \cite{GayBalmaz2017b} to include resistivity. This formulation involves three additional variables besides $ \varphi (t) \in \operatorname{Diff}_0( \Omega ) $, $ \varrho (t) {\rm d} x \in \Lambda ^3( \Omega )$, and $S(t){\rm d} x \in \Lambda ^3( \Omega )$: the internal entropy density variable $ \Sigma (t) {\rm d} x \in \Lambda ^3( \Omega )$, whose time rate of change is the internal entropy production, the thermal displacement $ \Gamma (t) \in \Lambda ^0 ( \Omega )$, whose time rate of change is the temperature, and the magnetic displacement $ \mathcal{Z} (t) \cdot {\rm d} X \in \Lambda ^1( \Omega )$, whose time rate of change is the magnetic field.

The variational principle reads as follows. Find the curves $\varphi: [t_0,t_1] \rightarrow  \operatorname{Diff}_0( \Omega )$, $S {\rm d} x, \Sigma {\rm d}  x:[t_0,t_1]\rightarrow  \Lambda ^3   ( \Omega)$, $\Gamma: [t_0,t_1]\rightarrow \mathcal{F} ( \Omega )$, $ \mathcal{B} \cdot {\rm d} A: [t_0,t_1] \rightarrow \Lambda  ^2 ( \Omega )$, and $ \mathcal{Z} \cdot {\rm d} X :[t_0,t_1] \rightarrow \Lambda ^1( \Omega )$ which are critical for the \textit{variational condition}
\begin{equation}
\label{eq:VP_fluid} 
\delta \int_{t_0}^{t_1} \Big[L\big(\varphi, \dot \varphi,  S, \mathcal{B} , \varrho _0\big) +  \int_ \Omega \big(\mathcal{B} \cdot \dot{\mathcal{Z}}+ (S- \Sigma ) \dot \Gamma  \big)\,{\rm d} X \Big]{\rm d}t =0
\end{equation}
subject to the \textit{phenomenological constraint}
\begin{equation}
\label{eq:KC_fluid}
\frac{\delta  L }{\delta  S}\dot \Sigma= - P: \nabla\dot \varphi  + J_S \cdot \nabla \dot\Gamma  +J_ \mathcal{B}  \cdot \operatorname{curl}\dot{ \mathcal{Z}}
\end{equation}  
and for variations subject to the \textit{variational constraint}
\begin{equation}
\label{eq:VC_fluid}
\frac{\delta L}{\delta  S}\delta \Sigma = - P : \nabla \delta  \varphi + J_S \cdot \nabla \delta \Gamma +  J_\mathcal{B}  \cdot \operatorname{curl} \delta \mathcal{Z}
\end{equation} 
with $\delta \varphi|_{t=t_0,t_1}= \delta\Gamma|_{t=t_0,t_1}=\delta \mathcal{Z} |_{t=t_0,t_1}=0$, and $\delta\varphi|_{\partial\Omega}=0$. In \eqref{eq:KC_fluid} and \eqref{eq:VC_fluid} $\frac{\delta L}{\delta S}$ is the functional derivative of $L$ with respect to $S$, and is identified with minus the temperature of the fluid, denoted $ \mathfrak{T} = - \frac{\delta L}{\delta S}$ in the Lagrangian description. Applying the variational formulation \eqref{eq:VP_fluid}--\eqref{eq:VC_fluid} yields the equations of motion in the Lagrangian description, along with the conditions
\[
\dot \Gamma = - \frac{\delta L}{\delta S}=\mathfrak{T}\quad \text{and}\quad \dot{\mathcal{Z}}= - \frac{\delta L}{\delta \mathcal{B}}=\mathfrak{H},
\]
which assign to $\Gamma$ and $\mathcal{Z}$ their physical interpretations as thermal and magnetic displacements, respectively.

\paragraph{Eulerian description.} We write $ \sigma $, $j_s$, and $j_B$ the Eulerian fluxes associated to $P$, $J_S$, and $J_ \mathcal{B} $ above. In particular, we have
\begin{equation}
j_B \cdot {\rm d} x= \varphi _* ( J_ \mathcal{B} \cdot {\rm d} X),
\end{equation}
see \cite{GayBalmaz2017b,Gawlik2024} for the other fluxes. Besides the relations \eqref{eq:push_forward}, the Eulerian variables associated to $ \Gamma $, $ \Sigma $, and $ \mathcal{Z} $ are also needed
\begin{equation}\label{eq:push_forward_2} 
\gamma = \varphi _* \Gamma , \quad \varsigma {\rm d} x= \varphi _*( \Sigma  {\rm d} X) , \quad Z \cdot {\rm d} x= \varphi  _*(\mathcal{Z} \cdot {\rm d} X).
\end{equation}
The Eulerian version of the principle \eqref{eq:VP_fluid}--\eqref{eq:VC_fluid} reads as follows. Find the curves $u: [t_0,t_1] \rightarrow  \mathfrak{X} _0( \Omega )$, $ \rho  \,{\rm d} x, s{\rm d} x, \varsigma {\rm d} x :[t_0,t_1] \rightarrow  \Lambda ^3 ( \Omega)$, $ \gamma : [t_0,t_1] \rightarrow \Lambda ^0( \Omega )$,  $B \cdot {\rm d} a:[t_0,t_1] \rightarrow \Lambda ^2 ( \Omega )$, and $ Z \cdot {\rm d} x:[t_0,t_1] \rightarrow \Lambda ^1 ( \Omega  )$ which are critical for the \textit{variational condition}
\begin{equation}\label{eq:VP_NSF_spatial}
\delta \int_{t_0}^{t_1}\Big[ \ell(u , \rho , s,B)+\int_ \Omega \big(B \cdot \mathfrak{D}_t Z + (s- \varsigma )D_t \gamma\big) \, {\rm d} x \Big] {\rm d}t=0,
\end{equation}
with the \textit{phenomenological constraint} and \textit{variational constraint} given by
\begin{subequations}
\label{eq:NSF_spatial}
\begin{align}
\label{eq:KC_NSF_spatial}
\frac{\delta   \ell}{\delta   s} \overline{D}_t \varsigma & =  -  \sigma : \nabla u + j_s \cdot \nabla (D_t \gamma) +j_B \cdot \operatorname{curl}(\mathfrak{D}_t Z), \\
\label{eq:VC_NSF_spatial}
\frac{\delta   \ell}{ \delta  s}  \overline{D}_\delta  \varsigma & =   -\sigma : \nabla  v +  j_s \cdot \nabla (D_ \delta \gamma)+ j_B \cdot \operatorname{curl}(\mathfrak{D}_\delta  Z),
\end{align}
\end{subequations}
and the Euler-Poincar\'e constraints
\begin{equation}\label{eq:EP_constraints} 
\delta u = \partial _t v - \pounds_v u, \qquad \overline{D}_\delta \rho = 0,
\end{equation}
with $v:[t_0,t_1] \rightarrow \mathfrak{X} ( \Omega )$ an arbitrary time-dependent vector field vanishing at the boundary. It eventually yields the system of equations
\begin{subequations}
\label{eq:eulerian} 
\begin{align}
\label{eq:eulerian_u} 
& \partial _t \frac{\delta  \ell}{\delta u} + \pounds _u \frac{\delta \ell}{\delta u} = \rho  \nabla \frac{\delta \ell}{\delta \rho  }+ s \nabla \frac{\delta \ell}{\delta s}+ \operatorname{div} \sigma +  B \times \operatorname{curl} \frac{\delta \ell}{\delta B}, \\
& \overline{\mathfrak{D}}_t B = \operatorname{curl}j_B, \\
& \overline{D}_t \rho  = 0, \\
\label{eq:eulerian_s} 
& \displaystyle \frac{\delta \ell}{\delta s} (\overline{D}_t s + \operatorname{div}j_s)= - \sigma : \nabla u - j_s \cdot  \nabla \frac{\delta \ell}{\delta s} - j_B \cdot \operatorname{curl} \frac{\delta  \ell}{\delta  B},
\end{align}
\end{subequations}
together with the boundary conditions
\begin{equation}\label{eq:BC} 
 j_s \cdot n=0, \quad j_B \times n=0,
\end{equation} 
and the conditions
\begin{equation}\label{eq:additional_cond} 
\overline{D}_t \varsigma = \overline{D}_t s + \operatorname{div}j_s, \quad D_t \gamma = - \frac{\delta \ell}{\delta s}, \quad \mathfrak{D}_t Z= - \frac{\delta \ell}{\delta B}.
\end{equation} 
The boundary conditions $u|_{ \partial \Omega }=0$ and $B \cdot n=0$ are assumed from the start. The derivation of \eqref{eq:eulerian}--\eqref{eq:additional_cond} from \eqref{eq:VP_NSF_spatial}--\eqref{eq:EP_constraints} follows a similar approach to that in \cite{GayBalmaz2017b} for the Navier-Stokes-Fourier case, to which we refer for further details.

\paragraph{Energy and entropy balance.} By defining the total energy $ \mathcal{E} = \left\langle \frac{\delta \ell}{\delta u}, u \right\rangle - \ell(u, \rho  , s, B)$, from \eqref{eq:eulerian} we get the energy conservation
\begin{equation}\label{eq:energy_cons} 
\frac{d}{dt} \mathcal{E} = \int_ \Omega \operatorname{div}\Big(\Big( - \frac{\delta \ell}{\delta u} \cdot u + \rho  \frac{\delta \ell}{\delta \rho  }  + s  \frac{\delta \ell}{\delta s  }\Big)u + \frac{\delta \ell}{\delta B} \times (u \times B) + \sigma \cdot u + j_s \frac{\delta \ell}{\delta s} + \frac{\delta \ell}{\delta B}  \times j_B \Big) {\rm d} x=0,
\end{equation}
where the last equality follows from the boundary conditions $u|_{ \partial \Omega }=0$, $j_s \cdot n=0$, and $j_B \times n=0$. As for the entropy equation, it reads
\begin{equation}
\label{eq:generic_entropy_prod}
\overline{D}_t s + \operatorname{div}j_s = \frac{1}{\frac{\delta \ell}{\delta s}} \left( - \sigma : \nabla u - j_s \cdot  \nabla \frac{\delta \ell}{\delta s} - j_B \cdot \operatorname{curl} \frac{\delta  \ell}{\delta  B} \right).
\end{equation}
Proper closures for $\sigma, j_s$ and $j_B$ must then be given so as to satisfy the Second Law of Thermodynamics (namely for the right-hand side of \eqref{eq:generic_entropy_prod} to be positive).

\paragraph{Treatment of other thermal boundary conditions.} The variational formulation given above yields thermally insulated boundaries, as we see from the first equation in \eqref{eq:BC}, consistently with the variational setting of adiabatically closed systems, \cite{GayBalmaz2017a}. As shown in \cite{Gawlik2024}, a modification of the constraints \eqref{eq:KC_NSF_spatial} and \eqref{eq:VC_NSF_spatial}, consistent with the variational principle for open systems, \cite{GayBalmaz2018b}, allows the treatment of the Dirichlet boundary condition $T|_{ \partial \Omega }=T_0$ or the nonhomogeneous Neumann boundary condition $j_q \cdot n=q_0$, with $j_q=Tj_s$ the heat flux, for given $T_0, q_0: \partial \Omega \rightarrow \mathbb{R} $. To obtain these boundary conditions, the constraint \eqref{eq:KC_NSF_spatial} is modified with the addition of a boundary term, namely,
\begin{equation}
\label{eq:KC_NSF_spatial_w_Dirichlet}
\int_ \Omega w\frac{\delta \ell}{\delta   s} \overline{D}_t \varsigma  \,{\rm d} x  = \int_ \Omega w\big(- \sigma : \nabla u +j_s \cdot \nabla D_t \gamma +  j_B \cdot \operatorname{curl}(\mathfrak{D}_t Z)\big)\, {\rm d} x - \int_ { \partial \Omega } \!\!\!w(j_s \cdot n) \left( D_ t \gamma  -T_0 \right)\,  {\rm d} a,\;\forall\; w
\end{equation}
for the Dirichlet boundary conditions, and
\begin{equation}
\label{eq:KC_NSF_spatial_w_nhom}
\int_ \Omega w\frac{\delta \ell}{\delta   s} \overline{D}_t \varsigma  \,{\rm d} x  = \int_ \Omega w\big(- \sigma : \nabla u +j_s \cdot \nabla D_t \gamma +  j_B \cdot \operatorname{curl}(\mathfrak{D}_t Z)\big)\, {\rm d} x - \int_ { \partial \Omega } w\big( (j_s \cdot n) D_ t \gamma  - q_0 \big)\, {\rm d} a,\;\forall\; w
\end{equation} 
for the nonhomogeneous Neumann boundary condition. For both cases, the variational constraint \eqref{eq:VC_NSF_spatial} is modified as
\begin{equation}\label{eq:VC_NSF_spatial_w_Dirichlet}
\int_ \Omega w\frac{\delta   \ell}{ \delta  s}  \overline{D}_\delta  \varsigma\, {\rm d} x  =\int_ \Omega w \big(- \sigma : \nabla  v + j_s \cdot \nabla D_ \delta \gamma + j_B \cdot \operatorname{curl}( \mathfrak{D}_\delta  Z)\big )\, {\rm d} x   - \int_ { \partial \Omega } w(j_s \cdot n) D_ \delta \gamma \,  {\rm d} a,\;\forall\;w.
\end{equation}
The variational principle \eqref{eq:VP_NSF_spatial}--\eqref{eq:EP_constraints} with \eqref{eq:KC_NSF_spatial} replaced by \eqref{eq:KC_NSF_spatial_w_Dirichlet} or \eqref{eq:KC_NSF_spatial_w_nhom}, and with \eqref{eq:VC_NSF_spatial} replaced by \eqref{eq:VC_NSF_spatial_w_Dirichlet}, yield the same system \eqref{eq:eulerian}--\eqref{eq:additional_cond}, but with the boundary condition $j_s \cdot n=0$ replaced by $T|_{ \partial \Omega }=T_0$ or $(j_s \cdot n) T =q_0$. In these cases, the energy balance in \eqref{eq:energy_cons} is modified as
\begin{subequations}\label{energy_balance}
\begin{align}
\frac{d}{dt} \mathcal{E} & = \int_{ \partial \Omega } \frac{\delta \ell}{\delta s} (j_s \cdot n) {\rm d} a= -  \int_{ \partial \Omega } T_0( j_s \cdot n ){\rm d} a, \\
\frac{d}{dt} \mathcal{E} & = \int_{ \partial \Omega } \frac{\delta \ell}{\delta s} (j_s \cdot n) {\rm d} a= -\int_{ \partial \Omega }q_0 {\rm d} a.
\end{align} 
\end{subequations}

\subsection{Physical closures}

\paragraph{MHD Lagrangian.} We consider the general case of MHD with Coriolis and gravitational forces, where the Lagrangian is given by
\begin{equation}
\label{eq:MHD_Lagrangian} 
\ell(u, \rho  , s, B)= \int_ \Omega \Big[\frac{1}{2} \rho  | u| ^2 + \rho  R \cdot u - \varepsilon ( \rho  , s) - \frac{1}{2 \mu_0} | B | ^2 - \rho \, \phi  \Big] {\rm d} x
\end{equation} 
with $R$ the vector potential of the angular velocity, $ \varepsilon ( \rho  , s)$ the internal energy density and $ \mu_0 >0$ the magnetic permeability. One then gets $\frac{\delta \ell}{\delta u}= \rho ( u+R)$, $\frac{\delta  \ell}{\delta  \rho  } = \frac{1}{2}| u | ^2 + R \cdot u - \frac{\partial \varepsilon }{\partial \rho  } - \phi $,  $-\frac{\delta \ell}{\delta s}= \frac{\partial \varepsilon  }{\partial s}= T >0$ the temperature, and $ -\frac{\delta \ell}{\delta B}= \frac{1}{ \mu _0}B=H$. By defining the pressure $p= \frac{\partial \varepsilon }{\partial \rho  } \rho  + \frac{\partial \varepsilon }{\partial s} s - \varepsilon$, \eqref{eq:eulerian_u} and \eqref{eq:eulerian_s} eventually become
\begin{subequations}
\begin{align}
\rho  ( \partial _t u+ u \cdot \nabla u +2 \omega \times u) & = - \rho  \nabla \phi - \nabla p + \operatorname{curl} H \times B+ \operatorname{div} \sigma, \\
T (\overline{D}_ts  + \operatorname{div} j _s) & =  \sigma \!: \!\nabla u - j _s\! \cdot  \!\nabla T -j_B \! \cdot \! \operatorname{curl}H.
\end{align}
\end{subequations}

\paragraph{Dissipation and entropy production.} We shall first consider the usual expressions
\begin{subequations}\label{eq:phenom_relations}
\begin{align} 
\sigma & = \sigma (u) = 2 \mu \operatorname{Def} u +  \lambda (\operatorname{div}u) \delta, & & \text{with $\mu \geq 0$ and $\zeta = \lambda + \frac{2}{d} \mu \geq 0$,}\\
\label{eq:phenom_relations_jb}
j_s & = j_s (T) =-  \frac{1}{T} \kappa \nabla T, & & \text{with $\kappa \geq 0$},\\
\label{eq:phenom_relations_js}
j_B & = j_B(H) =- \nu \operatorname{curl}H, & & \text{with $ \nu \geq 0$},
\end{align}
\end{subequations} 
where $\operatorname{Def}u = \frac{1}{2}(\nabla u + \nabla u^\mathsf{T}) $ is the rate of deformation tensor, $\mu\geq 0$ and $ \zeta = \lambda + \frac{2}{d} \mu \geq 0$ are the shear and bulk viscosity coefficients, $ \kappa\geq 0$ is the thermal conductivity coefficient, and $\nu \geq 0$ is the resistivity. In this paper, these coefficients will be assumed constant. Their sign is imposed by the Second Law of Thermodynamics $ \overline{D}_t s+ \operatorname{div}j_s \geq 0$. Indeed, with \eqref{eq:phenom_relations}, the entropy equation reads
\begin{align}
\notag
\overline{D}_t s+ \operatorname{div}j_s & =  \frac{1}{T} \sigma \!: \!\nabla u - \frac{1}{T} j _s\! \cdot  \!\nabla T - \frac{1}{ T}j_B \! \cdot \!   \operatorname{curl}H\\
\label{eq:entrop_prod} 
&= \frac{2 \mu }{T}\left( \operatorname{Def} u \right) ^{(0)} \!: \! \left( \operatorname{Def} u \right) ^{(0)} + \frac{\zeta}{T} (\operatorname{div}u) ^2  +  \frac{\kappa}{T ^2 } | \nabla T| ^2 + \frac{1}{T} \nu  | \operatorname{curl}H| ^2  \geq 0,
\end{align}
where $\left( \operatorname{Def} u \right) ^{(0)}$ denotes the trace-free part of $\operatorname{Def} u$.

\paragraph{Thermoelectric effect (Seebeck--Peltier coupling).}
It is possible to extend \eqref{eq:phenom_relations_jb}-- \eqref{eq:phenom_relations_js} to include the thermoelectric effect (Seebeck--Peltier coupling). The new closure reads
\begin{subequations}
\label{eq:phenom_relations_SP}
\begin{align} 
\label{eq:phenom_relations_jb_SP}
j_s & = j_s (T, H) =-  \frac{1}{T} \kappa \nabla T + \alpha \operatorname{curl} H, & & \text{with $\kappa \geq 0$}, \\
\label{eq:phenom_relations_js_SP}
j_B & = j_B(T, H) = - \alpha \nabla T - \nu \operatorname{curl}H, & & \text{with $ \nu \geq 0$},
\end{align}
\end{subequations} 
for some coefficient $\alpha$. Contrary to the previous coefficients, $\alpha$ is not subject to any sign restriction, as the thermoelectric effect only contributes to the entropy balance through a flux, the so-called Thomson effect. Note that if $\alpha$ is constant or a sole function of the temperature, it does not contribute to the induction equation.

\subsection{Associated weak formulation}

The variational derivation presented above yields a weak formulation where the evolution equations and boundary conditions are treated simultaneously. This formulation will be shown to have a discrete version, which allows to achieve thermodynamic consistency at the discrete level. 

\paragraph{Phenomenological and variational constraints.} The treatment of viscosity and heat conduction is taken from \cite{Gawlik2024}. For viscosity, the trilinear form $c : L^\infty(\Omega) \times H^1(\Omega)^3 \times H^1(\Omega)^3 \rightarrow \mathbb{R}$ is defined by
\begin{equation}\label{def_c} 
c(w,u,v) = \int_\Omega w \, \sigma (u) : \nabla v \, {\rm d}x.
\end{equation}
For heat conduction, we set $F = \{f \in H^1(\Omega) \mid 1/f \in L^\infty(\Omega), \, \dv(\nabla f/f) \in L^2(\Omega)\}$ and define $d : W^{1,\infty}(\Omega) \times F \times H^1(\Omega) \rightarrow \mathbb{R}$ by
\begin{equation}
\label{d}
d(w,f,g) = \left\{
\begin{array}{ll}
\vspace{0.2cm}\displaystyle\int_\Omega w \, j_s(f) \cdot \nabla g \, {\rm d} x & \text{for homogeneous Neumann}\\
\displaystyle \int_\Omega w \, j_s(f) \cdot \nabla g \, {\rm d} x - \int_{\partial\Omega} w j_s(f) \cdot n g \, {\rm d}a & \parbox[c]{5cm}{for Dirichlet and \\ nonhomogeneous Neumann}
\end{array}
\right.
\end{equation}
and $e :  W^{1,\infty}(\Omega) \times F \rightarrow \mathbb{R}$ by
\begin{equation}
\label{e}
e(w,f) = \left\{
\begin{array}{ll}
\vspace{0.2cm}\displaystyle 0 & \text{for homogeneous Neumann}\\
\vspace{-0.4cm}\displaystyle \int_{\partial\Omega} w j_s(f) \cdot n T_0 \, {\rm d}a & \text{for Dirichlet}\\
\\
\displaystyle  \int_{\partial\Omega} w q_0(f) \, {\rm d}a &  \text{for nonhomogeneous Neumann}.
\end{array}
\right.
\end{equation}
Finally, resistivity is included through the trilinear form
\begin{equation}\label{eq:g}
g(w,C,D)=\int_ \Omega w j_B(C) \cdot \operatorname{curl}D {\rm d} x .
\end{equation}
With these definitions, the constraints \eqref{eq:KC_NSF_spatial}--\eqref{eq:VC_NSF_spatial}, \eqref{eq:KC_NSF_spatial_w_Dirichlet}--\eqref{eq:VC_NSF_spatial_w_Dirichlet}, and \eqref{eq:KC_NSF_spatial_w_nhom}--\eqref{eq:VC_NSF_spatial_w_Dirichlet} can be written in a unified way for all boundary conditions as
\begin{subequations}
\begin{align}
\label{eq:unified_writing_KC} 
\Big\langle  w,\frac{\delta   \ell}{\delta   s} \overline{D}_t \varsigma \Big\rangle & =- c(w,u,u) + g\Big(w, - \frac{\delta \ell}{\delta B}, \overline{\mathfrak{D}}_tZ \Big) + d\Big( w, - \frac{\delta \ell}{\delta s}, D_t \gamma \Big)+ e\Big(w, - \frac{\delta \ell}{\delta s} \Big), & \forall w, \\
\label{eq:unified_writing_VC}  
\Big\langle w, \frac{\delta   \ell}{ \delta  s}  \overline{D}_\delta  \varsigma \Big\rangle & =- c(w,u, v ) + g\Big(w, - \frac{\delta \ell}{\delta B}, \overline{\mathfrak{D}}_\delta Z \Big) + d\Big( w, - \frac{\delta \ell}{\delta s}, D_\delta  \gamma \Big), & \forall w,
\end{align}
\end{subequations} 
with $ \left\langle \cdot , \cdot \right\rangle $ the $L^2$ inner product.

\paragraph{Weak form of the variational principle.} By using these notations when carrying out the variational formulations \eqref{eq:VP_NSF_spatial}--\eqref{eq:EP_constraints} or \eqref{eq:VP_NSF_spatial}-\eqref{eq:KC_NSF_spatial_w_Dirichlet}-\eqref{eq:VC_NSF_spatial_w_Dirichlet}-\eqref{eq:EP_constraints} or \eqref{eq:VP_NSF_spatial}-\eqref{eq:KC_NSF_spatial_w_nhom}-\eqref{eq:VC_NSF_spatial_w_Dirichlet}-\eqref{eq:EP_constraints}, we get the following weak form of the equations \eqref{eq:eulerian} and associated boundary conditions $j_s \cdot n=0$ or $T|_{ \partial \Omega }=T_0$ or $j_q \cdot n=q_0$ as follows
\begin{subequations}
\label{eq:good_weak_form}
\begin{align}
\label{eq:good_weak_form_u}
& \left\langle \partial _t \frac{\delta \ell}{\delta u}, v  \right\rangle + a \Big( \frac{\delta \ell}{\delta u} , u, v \Big) + b ^1 \Big( \frac{\delta \ell}{\delta \rho  } , \rho, v \Big) +b ^1 \Big( \frac{\delta \ell}{\delta s  } , s  , v \Big) + b ^2 \Big( \frac{\delta \ell}{\delta B  } , B  , v \Big) = - c(1,u,v), & \forall v, \\
\notag
& - \left\langle \partial _t  s, \frac{\delta \ell}{\delta s}w \right\rangle + b ^1  \Big( -\frac{\delta \ell}{\delta s}w, s, u\Big) - d \Big( 1, -\frac{\delta \ell}{\delta s}, - \frac{\delta \ell}{\delta s} w\Big)
= c(w, u,u ) & \\
\label{eq:good_weak_form_s}
& \hspace{5cm} - g\Big( w, -\frac{\delta \ell}{\delta B} , - \frac{\delta \ell}{\delta B} \Big) - d\Big( w, -\frac{\delta \ell}{\delta s}, - \frac{\delta \ell}{\delta s}\Big)- e\Big(w, - \frac{\delta \ell}{\delta s} \Big), & \forall w,  \\
\label{eq:good_weak_form_rho}
& \left\langle \partial _t \rho  , \theta  \right\rangle + b^1(\theta , \rho  , u)=0, & \forall \theta, \\
\label{eq:good_weak_form_B}
& \left\langle \partial _t B, C \right\rangle +b^2(C,B,u)= g\Big(1,- \frac{\delta \ell}{\delta B} , C\Big), & \forall C.
\end{align}
\end{subequations}
The thermal boundary condition is weakly enforced by the entropy equation \eqref{eq:good_weak_form_s}, while the boundary condition $j_B \times n=0$ is a consequence of the weak version of the induction equation \eqref{eq:good_weak_form_B}.

\section{Structure preserving discretization}\label{sec_3}

We will now construct a spatial discretization of \eqref{eq:good_weak_form} using finite elements. It is based on a discrete variational principle and preserves both the energy conservation and the positivity of internal entropy production at the discrete level.

\subsection{Discrete setting}

We make use of the following function spaces
\begin{subequations}
\begin{align}
H^1_0(\Omega) & = \{ f \in L^2 (\Omega) \, | \, \nabla f \in L^2(\Omega)^3, \, f = 0 \text{ on } \partial \Omega \}, \\
H_0 (\operatorname{div}, \Omega) & = \{ u \in L^2 (\Omega)^3 \, | \operatorname{div} u \in L^2 (\Omega), \, u \cdot n = 0 \text{ on } \partial \Omega \}, \\
H_0 (\operatorname{curl}, \Omega) & = \{ u \in L^2 (\Omega)^3 \, | \, \operatorname{curl} u \in L^2(\Omega)^3, \, u \times n = 0 \text{ on } \partial \Omega \}.
\end{align}
\end{subequations}
Let $\mathcal{T}_h$ be a triangulation of $\Omega$. We regard $\mathcal{T}_h$ as a member of a family of triangulations parametrized by $h = \operatorname{max}_{K \in \mathcal{T}_h} h_K$, where $h_K$ denotes the diameter of a simplex $K$. We assume that this family is shape-regular, meaning that the ratio $\operatorname{max}_{K \in \mathcal{T}_h} h_K /\rho_K$ is bounded above by a positive constant for all $h > 0$. Here, $\rho_K$ denotes the inradius of $K$. 

When $r \geq 0$ is an integer and $K$ is a simplex, we write $P_r(K)$ to denote the space of polynomials on $K$ of degree at most $r$. Let $r, s \geq 0$ be fixed integers. The velocity $u$ and magnetic field $B$ are respectively discretized by continuous Galerkin and Raviart-Thomas finite elements. The density $\rho$ and entropy $s$ are both discretized with the discontinuous Galerkin spaces.
\begin{subequations}
\begin{align}
U^\text{grad}_h = CG_{r+1} (\mathcal{T}_h)^3 & \vcentcolon = \{ u \in H^1_0 (\Omega)^3 \, | \, u|_K \in P_{r+1}(K)^3, \, \forall K \in \mathcal{T}_h \} \\ 
U^\text{div}_h = RT_{r} (\mathcal{T}_h) & \vcentcolon = \{ B \in H_0(\operatorname{div}, \Omega) \, | \, B|_K \in P_{r}(K)^3 + xP_{r}(K), \, \forall K \in \mathcal{T}_h \} \\ 
F_h = DG_{s} (\mathcal{T}_h) & \vcentcolon = \{ f \in L^2(\Omega) \, | \, f|_K \in P_{s}(K), \, \forall K \in \mathcal{T}_h \}.
\end{align}
\end{subequations}
An auxiliary space will also be needed, namely the Nedelec finite element space of the first kind
\begin{equation}
U^\text{curl}_h = NED_r(\mathcal{T}_h) \vcentcolon = \{ u \in H_0(\operatorname{curl}, \Omega) \, | \, u|_K \in P_{r}(K)^3 + x \times P_r(K)^3, \, \forall K \in \mathcal{T}_h \}
\end{equation}
which satisfies $\operatorname{curl} U^{\operatorname{curl}}_h \subset U^\text{div}_h$. The $L^2(\Omega)$-orthogonal projections on $F_h$, $U^{\operatorname{div}}_h$, $U^{\operatorname{grad}}_h$, and $U^{\operatorname{curl}}_h$ are respectively denoted with $\pi_h$, $\pi^{\operatorname{div}}_h$, $\pi^{\operatorname{grad}}_h$, and $\pi^{\operatorname{curl}}_h$. Because discontinuous elements are used, some notations need to be introduced. Let $\mathcal{E}_h = \mathcal{E}^0_h \cup \mathcal{E}^\delta_h$ denote the set of codimension-$1$ faces in $\mathcal{T}_h$ with $\mathcal{E}^0_h$ the set of interior faces and $\mathcal{E}^\delta_h$ the set of boundary faces. For every face $e = K_1 \cap K_2 \in \mathcal{E}^0_h$, its length is written $h_e$, while the jump and average of a piecewise smooth scalar function $f$ are defined as
\begin{subequations}
\begin{align}
\llbracket f \rrbracket & = f_1 n_1 + f_2 n_2, \\
\{ f \} & = \frac{1}{2} (f_1 + f_2),
\end{align}
\end{subequations}
where $f_i = f |_{K_i}$ and $n_i$ is the normal vector to $e$ pointing outward from $K_i$. 

\begin{remark}
For 2D computations, the space $H_0(\operatorname{curl}, \Omega)$ and $U^{\text{curl}}_h$ respectively become
\begin{subequations}
\begin{align}
H_0 (\operatorname{curl}, \Omega) & = \{ u \in L^2 (\Omega)^2 \, | \, \partial_x u_y - \partial_y u_x \in L^2(\Omega), \, u_x n_y - u_y n_x = 0 \text{ on } \partial \Omega \}, \\
\label{eq:Ucurl_2D}
U^\text{curl}_h & = \{ u \in H_0(\operatorname{curl}, \Omega) \, | \, u|_K \in P_{r}(K)^2 + (y, -x) P_r(K), \, \forall K \in \mathcal{T}_h \}.
\end{align}
\end{subequations}
\end{remark}

\subsection{Discrete variational formulation}

The case without thermoelectric effect is first considered.

\paragraph{Discrete advection operators.}
Before giving the discrete weak formulation, it is necessary to define discrete counterparts of the advection operators $a$, $b^1$ and $b^2$. They are here chosen as
\begin{subequations}
\begin{align}
a_h(w,u,v) & = a(w,u,v) = -\int_\Omega w \cdot [u, v] {\rm d} x \\
b^1_h(f,g,u) & = - \sum_{K \in \mathcal{T}_h} \int_K (u \cdot \nabla f)g \text{d} x + \sum_{e \in \mathcal{E}_h} \int_e u \cdot \llbracket f \rrbracket \{ g \} {\rm d} s \\
b^2_h(C,B,u) & = \left \langle C, \operatorname{curl} \pi^{\operatorname{curl}}_h \left( \pi^{\operatorname{curl}}_h B \times u \right) \right \rangle.
\end{align}
\end{subequations}
Here $b^1_h$ is a standard discretization \cite{Brezzi2004}, while the choice of $b^2_h$ is motivated by conservation of magnetic helicity \cite{Gawlik2021}. Note that in the definition of $b^1_h$, a centered flux $\{ g \}$ was chosen in the boundary term for simplicity. Other choices are also possible; see, for example, \cite{Gawlik2024}[\S 3.5] for a description of how to instead use an upwinded flux without compromising structure preservation.

\paragraph{Discrete dissipation operators.} Discrete versions of operators $c, d, e$ and $g$ are chosen consistently with \cite{Gawlik2021} and \cite{Gawlik2024}. The viscosity operator $c_h : F_h \times U^{\operatorname{grad}}_h \times U^{\operatorname{grad}}_h \rightarrow \mathbb{R}$ and resistivity operator $g_h : F_h \times U^{\operatorname{div}}_h \times U^{\operatorname{div}}_h$  are simply given as
\begin{subequations}
\begin{align}
c_h(w,u,v) & = c(w,u,v) = \int_\Omega w \sigma(u) : \nabla v {\rm d} x, \\
\label{eq:g_h}
g_h(w,C,D) & = \int_\Omega \nu (\pi^{\operatorname{curl}}_h \operatorname{curl} C) \cdot \operatorname{curl} D {\rm d} x,
\end{align}
\end{subequations}
where the curl operator is interpreted in a distributional sense in \eqref{eq:g_h}. The operators $d_h : F_h \times F_h \times F_h\rightarrow\mathbb{R}$ and $e_h: F_h \times F_h \rightarrow\mathbb{R}$ both depend on the prescribed boundary conditions. For homogeneous Neumann conditions, they read
\begin{subequations}
\begin{align}
\notag
d_h^N(w, f, g) & = - \sum_{K \in \mathcal{T}_h} \int_K \frac{w}{f} \kappa \nabla f \cdot \nabla g {\rm d} x + \sum_{e \in \mathcal{E}_h^0} \int_e \frac{1}{\{f\}} \{ w \kappa \nabla f \} \cdot \llbracket g \rrbracket {\rm d} a \\
& {\color{white} = } - \sum_{e \in \mathcal{E}_h^0} \int_e \frac{1}{\{f\}} \{ w \kappa \nabla g \} \cdot \llbracket f \rrbracket {\rm d} a - \sum_{e \in \mathcal{E}_h^0} \frac{\eta}{h_e} \int_e \frac{\{w\}}{\{f\}} \llbracket f \rrbracket \cdot \llbracket g \rrbracket {\rm d} a, \\ 
e_h^{N}(w, f) & = 0,
\end{align}
\end{subequations}
where $\eta > 0$ is a penalty parameter. It is a standard non-symmetric interior penalty discretization of the Laplacian \cite[§ 10.5]{Brenner2008}. It will be shown, in Section \ref{sec:properties}, to be compatible with the Second Law of Thermodynamics. For non-homogeneous Neumann boundary conditions, they are defined as
\begin{subequations}
\begin{align}
d_h^{NN}(w, f, g) & = d_h^{N}(w, f, g) + \int_{\partial \Omega} \frac{w}{f} \kappa \nabla f \cdot n g {\rm d} a, \\ 
e_h^{NN}(w, f) & = \int_{\partial \Omega} w q_0 {\rm d} a.
\end{align}
\end{subequations}
Finally, for Dirichlet boundary conditions, the following discretization is made
\begin{subequations}
\begin{align}
d_h^{D}(w, f, g) & = d_h^{N}(w, f, g) - \int_{\partial \Omega} \frac{w}{f} \kappa \nabla g \cdot n (f - T_0) {\rm d} a + \int_{\partial \Omega} \frac{w}{f} \kappa \nabla f \cdot n g {\rm d} a,  \\ 
e_h^{D}(w, f) & = \int_{\partial \Omega} \frac{w}{f} \kappa \nabla f \cdot n T_0 {\rm d} a + \sum_{e \in \mathcal{E}^\partial_h} \frac{\eta}{h_e} \int_e w (f-T_0) {\rm d} a.
\end{align}
\end{subequations}

\paragraph{Numerical scheme.}
The discrete version of \eqref{eq:good_weak_form} then consists in finding $u \in U^{\operatorname{grad}}_h$, $B \in U^{\operatorname{div}}_h$, and $\rho, s \in F_h$ such that 
\begin{subequations}
\label{eq:scheme}
\begin{align}
& \left\langle \partial_t \frac{\delta \ell}{\delta u}, v \right\rangle + a\left( \pi^{\operatorname{grad}}_h \frac{\delta \ell}{\delta u}, u, v \right) + b^1_h\left( \pi_h \frac{\delta \ell}{\delta \rho}, \rho, v \right) + b^1_h\left( \pi_h \frac{\delta \ell}{\delta s}, s, v \right) + b^2_h\left( \pi^{\operatorname{div}}_h \frac{\delta \ell}{\delta B}, B, v \right) \\
\label{eq:scheme_u}
&\;\;\;= -c(1,u,v)\\
\label{eq:scheme_B}
& \left \langle \partial_t B, C \right \rangle + b^2_h\left( C, B, u \right) = g_h \left( 1, - \pi^{\operatorname{div}}_h \frac{\delta \ell}{\delta B}, C \right), \\
\label{eq:scheme_rho}
& \left \langle \partial_t \rho, \sigma \right \rangle + b^1_h\left( \sigma, \rho, u \right) = 0, \\
& - \left\langle \partial_t s, w \, \pi_h \frac{\delta \ell}{\delta s} \right\rangle + b^1_h \left( - \pi_h \left( w \, \pi_h \frac{\delta \ell}{\delta s} \right), s, u \right) - d_h \left(1, -\pi_h \frac{\delta \ell}{\delta s}, - \pi_h \left( w \, \pi_h \frac{\delta \ell}{\delta s} \right) \right) \notag \\ 
\label{eq:scheme_s}
&\;\;\; = c(w, u, u) - d_h \left(w, -\pi_h \frac{\delta \ell}{\delta s}, - \pi_h \frac{\delta \ell}{\delta s} \right) - g_h \left( w, - \pi^{\operatorname{div}}_h \frac{\delta \ell}{\delta B}, - \pi^{\operatorname{div}}_h \frac{\delta \ell}{\delta B} \right) - e_h\left( w, - \pi_h \frac{\delta \ell}{\delta s} \right),
\end{align}
\end{subequations}
for all $(v, C, \sigma, w) \in U^{\operatorname{grad}}_h \times U^{\operatorname{div}}_h \times F_h \times F_h$.

\begin{remark}
In the absence of a magnetic field, the numerical scheme can be derived from a discrete version of the variational principle \eqref{eq:VP_NSF_spatial}-\eqref{eq:NSF_spatial}-\eqref{eq:EP_constraints} as in \cite{Gawlik2024}. Members of $U^{\operatorname{grad}}_h$ are then seen as discrete vector fields that act on functions and densities (i.e. $0$-forms and $3$-forms), see \cite{GaGB2021}. In \cite{Carlier2025}, the action of vector fields is defined on forms of arbitrary order. It would allow to extend the discrete variational principle of \cite{Gawlik2024} to the full system of viscous-resistive MHD. The resulting scheme would only differ from \eqref{eq:scheme} in the expression of $b^2_h$ so as to reflect the discrete action of vector fields on $1$-forms and $2$-forms.
\end{remark}

In practice, it is possible to remove the projection of $\frac{\delta \ell}{\delta u}$ onto $U^{\operatorname{grad}}_h$ in \eqref{eq:scheme_u} without altering the properties of the scheme. Likewise for the outermost projection of $\pi_h \left( w \pi_h \frac{\delta \ell}{\delta s} \right)$ in \eqref{eq:scheme_s} and the projections $\pi^{\operatorname{div}}_h$ inside function $g_h$ in \eqref{eq:scheme_B} and \eqref{eq:scheme_s}. Because of the remaining projections, \eqref{eq:scheme} is not directly implementable in this form. In the spirit of \cite{Gawlik2021}, intermediate variables are added. For the Lagrangian \eqref{eq:MHD_Lagrangian}, \eqref{eq:scheme} is then equivalent to finding $u \in U^{\operatorname{grad}}_h$, $B \in U^{\operatorname{div}}_h$, $\rho, s, \theta, T \in F_h$, and $J, H, E \in U^{\operatorname{curl}}_h$ such that 
\begin{subequations}
\label{eq:scheme_final}
\begin{align}
\label{eq:scheme_final_u}
& \langle \partial_t (\rho u), v \rangle + a(\rho u, u, v) + b^1_h(\theta, \rho, v) + b^1_h(T, s, v) - \langle J \times H, v \rangle = -c(1, u, v), & \forall v \in U^{\operatorname{grad}}_h, \\
\label{eq:scheme_final_B}
& \langle \partial_t B, C \rangle + \langle \operatorname{curl} E, C \rangle = -\nu \langle \operatorname{curl} J, C \rangle & \forall C \in U^{\operatorname{div}}_h, \\
\label{eq:scheme_final_rho}
& \langle \partial_t \rho, \sigma \rangle + b^1_h(\sigma, \rho, u) = 0, & \forall \sigma \in F_h \\
\notag & \langle \partial_t s, T w \rangle + b^1_h(Tw, s, u) - d_h(1, T, Tw) \\ 
\label{eq:scheme_final_s}
& \hspace{0.15\textwidth} = c(w, u, u) + \nu w \langle J, J \rangle - d_h(w, T, T) - e_h(w, T), & \forall w \in F_h, \\
& \langle \theta, \tau \rangle = \left\langle \frac{\delta \ell}{\delta \rho}, \tau \right\rangle, & \forall \tau \in F_h, \\
& \langle T, \gamma \rangle = \left\langle - \frac{\delta \ell}{\delta s}, \gamma \right\rangle, & \forall \gamma \in F_h, \\
& \langle J, K \rangle = \frac{1}{\mu_0} \langle B, \operatorname{curl} K \rangle, & \forall K \in U^{\operatorname{curl}}_h, \\
& \langle H, G \rangle = \langle B, G \rangle, & \forall G \in U^{\operatorname{curl}}_h, \\
& \langle E, F \rangle = -\langle u \times H, F \rangle, & \forall F \in U^{\operatorname{curl}}_h.
\end{align}
\end{subequations}

\begin{remark}
For 2D computations, two cases are distinguished
\begin{itemize}
\item if the magnetic field is parallel to the plane, so is the variable $H$, while $J$ and $E$ are both orthogonal to it.
\item if the magnetic field is orthogonal to the plane, the opposite situation occurs: $H$ is orthogonal to the plane, while $J$ and $E$ are both parallel to it.
\end{itemize}
In both cases, parallel vectors are discretized with $U^{\operatorname{curl}}_h$ \eqref{eq:Ucurl_2D} while orthogonal vectors are discretized with 1D continuous Galerkin elements
\begin{equation}
CG_{r+1} (\mathcal{T}_h) = \{ u \in H^1_0 (\Omega) \, | \, u|_K \in P_{r+1}(K), \, \forall K \in \mathcal{T}_h \}.
\end{equation}
\end{remark}

\subsection{Discretization of the thermoelectric effect}

The thermoelectric coupling is incorporated by substituting the functions $d_h$ and $g_h$ with $\overline{d}_h : F_h \times F_h \times F_h \times U^{\operatorname{div}}_h \rightarrow \mathbb{R}$ and $\overline{g}_h : F_h \times U^{\operatorname{div}}_h \times U^{\operatorname{div}}_h \times F_h \rightarrow \mathbb{R}$
\begin{subequations}
\begin{align}
\overline{d}_h(w, f, g, D) & = d_h(w, f, g) + h_h(w, g, D), \\
\overline{g}_h(w, C, D, g) & = g_h(w, C, D) - h_h(w, g, D),
\end{align}
\end{subequations}
where $h_h : F_h \times F_h \times U^{\operatorname{div}}_h$ is a discretization of
\begin{subequations}
\begin{align}
h(w,g,D) & = \int_\Omega w \alpha \operatorname{curl} D \cdot \nabla g {\rm d} x \\
& = -\int_{\partial\Omega} D \cdot (w \alpha n \times \nabla g) \, {\rm d} a + \int_\Omega D \cdot (\nabla (w \alpha) \times \nabla g) \, {\rm d} x.
\end{align}
\end{subequations}
The following expression is here considered
\begin{multline}
h_h(w,g,D) = -\int_{\partial\Omega} D \cdot (w \alpha n \times \nabla g) \, {\rm d} a +  \sum_{K \in \mathcal{T}_h} \int_K D \cdot (\nabla (w\alpha) \times \nabla g) \, {\rm d} x \\
- \sum_{e \in \mathcal{E}_h^0} \int_e \{D\} \cdot (\llbracket w\alpha \rrbracket \times \{ \nabla g \}) \, {\rm d} a - \sum_{e \in \mathcal{E}_h^0} \int_e \{D\} \cdot (\{ \nabla(w\alpha) \} \times \llbracket g \rrbracket) \, {\rm d} a.
\end{multline}
The proof of the its consistency can be found in Appendix \ref{app:consistency}.

\subsection{Properties of the semi-discrete scheme}
\label{sec:properties}

The semi-discrete scheme \eqref{eq:scheme} reproduces at the discrete level several important features of continuous non-ideal MHD flows.

\paragraph{Mass conservation.}
Taking $\sigma = 1$ in the density equation \eqref{eq:scheme_rho} yields
\begin{equation}
\frac{\text{d}}{\text{d} t} \int_\Omega \rho \, {\rm d} x = \langle \partial_t \rho, 1 \rangle = -b_h^1(1, \rho, u) = 0.
\end{equation}

\paragraph{Gauss' magnetic law.}
The finite elements have been defined such that $\operatorname{curl} U^{\operatorname{curl}}_h \subset U^{\operatorname{div}}_h$ so that the following equation holds pointwise
\begin{equation}
\partial_t B + \operatorname{curl} E = - \nu \operatorname{curl} J.
\end{equation}
Then, following from the property $\operatorname{div} \circ \operatorname{curl} = 0$, it gives $\partial_t \operatorname{div} B = 0$. As such, if Gauss' magnetic law $\operatorname{div} B = 0$ is initially satisfied, it remains so during the entire computation.

\paragraph{Magnetic helicity balance.} Let $A$ be a vector satisfying $\operatorname{curl} A = B$ and $\left. A \times n \right|_{\partial \Omega} = 0$. Then the evolution of the helicity is given by
\begin{subequations}
\begin{align}
\frac{{\rm d}}{{\rm d}t} \langle A, B\rangle & = \langle \partial_t A, B\rangle + \langle \partial_t B, A \rangle \\
& = \langle \partial_t A, \operatorname{curl} A \rangle + \langle \partial_t B, A \rangle \\
& = \langle \operatorname{curl} \partial_t A, A \rangle + \langle \partial_t B, A \rangle \\
& = 2 \langle \partial_t B, A \rangle \\
& = - 2b^2_h(\pi^{\operatorname{div}}_h A, B, u) + 2g_h\left(1, -\pi^{\operatorname{div}}_h \frac{\delta \ell}{\delta B}, \pi^{\operatorname{div}}_h A\right).
\end{align}
\end{subequations}
The function $b^2_h$ has been defined to cancel in such a case. Indeed
\begin{subequations}
\begin{align}
b^2_h(\pi^{\operatorname{div}}_h A, B, u) & = \left \langle \pi^{\operatorname{div}}_h A, \operatorname{curl} \pi^{\operatorname{curl}}_h \left( \pi^{\operatorname{curl}}_h B \times u \right) \right \rangle \\
& = \left \langle A, \operatorname{curl} \pi^{\operatorname{curl}}_h \left( \pi^{\operatorname{curl}}_h B \times u \right) \right \rangle \\
& = \left\langle \pi^{\operatorname{curl}}_h \operatorname{curl} A,  \pi^{\operatorname{curl}}_h B \times u \right \rangle \\
& = \left\langle \pi^{\operatorname{curl}}_h B, \pi^{\operatorname{curl}}_h B \times u  \right \rangle = 0.
\end{align}
\end{subequations}
In particular, magnetic helicity is conserved in the absence of resistivity and thermoelectric coupling.

\paragraph{Energy balance.}
Likewise, taking $v = u$ in \eqref{eq:scheme_final_u}, $C = -\pi^{\operatorname{div}}_h \frac{\delta \ell}{\delta B}$ in \eqref{eq:scheme_final_B}, $\sigma = -\pi_h \frac{\delta \ell}{\delta \rho}$ in \eqref{eq:scheme_final_rho}, and $w = 1$ in \eqref{eq:scheme_final_s}, eventually gives the total energy rate of change
\begin{subequations}
\begin{align}
\frac{\text{d}}{\text{d} t} \int_\Omega \mathcal{E} \, {\rm d} x & = \left\langle \partial_t \frac{\delta \ell}{\delta u}, u \right\rangle - \left \langle \partial_t B, \pi^{\operatorname{div}}_h \frac{\delta \ell}{\delta B} \right \rangle - \left \langle \partial_t \rho, \pi_h \frac{\delta \ell}{\delta \rho} \right \rangle - \left \langle \partial_t s, \pi_h \frac{\delta \ell}{\delta s} \right \rangle \\
& = - e_h \left( 1, -\pi_h \frac{\delta \ell}{\delta s} \right).
\end{align}
\end{subequations}
In particular, total energy is conserved for homogeneous Neumann boundary conditions, and the energy balances in \eqref{energy_balance} are reproduced at the discrete level for the other boundary conditions.

\paragraph{Discrete Second Law of Thermodynamics.}
Let's assume homogeneous Neumann boundary conditions so that the system is insulated. Then consistency with thermodynamics requires that the right-hand side of \eqref{eq:scheme_final_s} be positive for all $w > 0$. This is indeed the case as
\begin{subequations}
\begin{align}
c(w, u, u) = \int_\Omega w \sigma(u) : u {\rm d} x & \geq 0, \\
- d^N_h(w, T, T) = \sum_{K \in \mathcal{T}_h} \frac{w}{T} |\nabla T|^2 {\rm d} x + \sum_{e \in \mathcal{E}^0_h} \frac{\eta}{h_e} \int_e \frac{\{w\}}{\{T\}} \left| \llbracket T \rrbracket \right|^2 {\rm d} a & \geq 0, \\
\nu w \langle J, J \rangle & \geq 0.
\end{align}
\end{subequations}
The same result still holds for non-homogeneous Neumann or Dirichlet boundary conditions as long as $w$ has compact support inside $\Omega$. It is also still valid when taking into account the thermoelectric effect, as it does not contribute to the right-hand side of the entropy equation.

\subsection{Time discretization}

The energy-conserving time discretization from \cite{Gawlik2024} is adapted to the present scheme. Solutions are approximated at discrete times $t_k$. The values at $t_k$ of quantity $\phi$ is $\phi_{k}$ while $\phi_{k+\frac{1}{2}}$ denotes the average $\frac{1}{2}(\phi_{k} + \phi_{k+1})$. The following discrete derivations of the energy will also be needed
\begin{subequations}
\begin{align}
\delta_\rho = \frac{1}{2} \left( \frac{\varepsilon(\rho_{k+1}, s_{k+1}) - \varepsilon(\rho_{k}, s_{k+1})}{\rho_{k+1} - \rho_k} + \frac{\varepsilon(\rho_{k+1}, s_k) - \varepsilon(\rho_{k}, s_k)}{\rho_{k+1} - \rho_k} \right), \\
\delta_s = \frac{1}{2} \left( \frac{\varepsilon(\rho_{k+1}, s_{k+1}) - \varepsilon(\rho_{k+1}, s_k)}{s_{k+1} - s_k} + \frac{\varepsilon(\rho_k, s_{k+1}) - \varepsilon(\rho_{k}, s_k)}{s_{k+1} - s_k} \right).
\end{align}
\end{subequations}
Given $u, B, \rho$ and $s$ at time $t_k$, their value at time $t_{k+1}$, with $t_{k+1} - t_k = \Delta t$, is obtained with the following fully discrete scheme
\begin{subequations}
\label{eq:time_discrete}
\begin{align}
\label{eq:time_discrete_u}
\notag
& \left\langle \frac{\rho_{k+1} u_{k+1} - \rho_k u_k}{\Delta t} , v \right\rangle + a((\rho u)_{k+\frac{1}{2}}, u_{k+\frac{1}{2}}, v) + b^1_h(\theta, \rho_{k+\frac{1}{2}}, v) & \\
& \hspace{0.1\textwidth} + b^1_h(T, s_{k+\frac{1}{2}}, v) - \langle J \times H, v \rangle = -c(1, u_{k+\frac{1}{2}}, v), & \forall v \in U^{\operatorname{grad}}_h, \\
\label{eq:time_discrete_B}
& \left\langle \frac{B_{k+1} - B_k}{\Delta t}, C \right\rangle + \langle \operatorname{curl} E, C \rangle = -\nu \langle \operatorname{curl} J, C \rangle - h_h(1, T, C), & \forall C \in U^{\operatorname{div}}_h, \\
\label{eq:time_discrete_rho}
& \left\langle \frac{\rho_{k+1} - \rho_k}{\Delta t}, \sigma \right\rangle + b^1_h(\sigma, \rho_{k+\frac{1}{2}}, u_{k+\frac{1}{2}}) = 0, & \forall \sigma \in F_h, \\
\notag & \left\langle \frac{s_{k+1} - s_k}{\Delta t}, T w \right\rangle + b^1_h(Tw, s, u_{k+\frac{1}{2}}) - d_h(1, T, Tw) - h_h\left(1, T \omega, \frac{1}{\mu_0} B_{k+\frac{1}{2}}\right) \\ 
\label{eq:time_discrete_s}
& \hspace{0.1\textwidth} = c(w, u_{k+\frac{1}{2}}, u_{k+\frac{1}{2}}) + \nu w \langle J, J \rangle - d_h(w, T, T) - e_h(w, T), & \forall w \in F_h, \\
& \langle \theta, \tau \rangle = \left\langle \frac{1}{2} u_k \cdot u_{k+1} - \delta_\rho, \tau \right\rangle, & \forall \tau \in F_h, \\
& \langle T, \gamma \rangle = \left\langle -\delta_s, \gamma \right\rangle, & \forall \gamma \in F_h, \\
\label{eq:time_discrete_J}
& \langle J, K \rangle = \frac{1}{\mu_0} \langle B_{k+\frac{1}{2}}, \operatorname{curl} K \rangle, & \forall K \in U^{\operatorname{curl}}_h, \\
& \langle H, G \rangle = \langle B_{k+\frac{1}{2}}, G \rangle, & \forall G \in U^{\operatorname{curl}}_h, \\
\label{eq:time_discrete_E}
& \langle E, F \rangle = -\langle u \times H, F \rangle, & \forall F \in U^{\operatorname{curl}}_h.
\end{align}
\end{subequations}
The proofs of Section \ref{sec:properties} can be adapted to show that the fully discrete scheme retains all the properties achieved so far, namely conservation of mass, Gauss' magnetic law, magnetic helicity balance, energy conservation, and consistency with the Second Law of Thermodynamics. Additionally, if dissipation is completely removed, the scheme is symmetric with respect to time inversion $t_k \leftrightarrow t_{k+1}$ showing that it is fully reversible.

\begin{remark}
Other time discretizations are also possible. Examples include the finite-element-in-time discretization of \cite{Andrews2025}, as well as an extension of the variational time discretization of \cite{GaGB2021}. The former uses auxiliary variables to systematically preserves invariants and would also allow to retain all properties of the present scheme.
\end{remark}

\begin{remark}
In our implementation, Newton's method is used to solve system \eqref{eq:time_discrete} at each time step.  Within each iteration of Newton's method, GMRES is used to solve the linear system.  A Schur complement preconditioner is also added, where the velocity is split from all other fields. The inverse of the Schur complement is approximated with a single multigrid V-cycle, and the inverse of the velocity block is approximated with an $ILU(0)$ preconditioner.
\end{remark}

\begin{remark}
Note that equations \eqref{eq:time_discrete_J}-\eqref{eq:time_discrete_E} form a linear systems in the unknowns $J$, $H$, and $E$ which involve non-diagonal mass matrices. To reduce computational expenses, it may be possible to replace these matrices with lumped mass matrices, effectively modifying the inner product on $U_h^{\mathrm{curl}}$. However, doing so alters the structure-preserving properties of the scheme; namely, conservation of helicity is lost in the reversible case, while a discrete second principle of thermodynamics is not as immediate in the irreversible case. An interesting topic of future research would be to design mass-lumping schemes that interfere minimally with structure preservation.
\end{remark}

\section{Numerical results}\label{sec_4}

For the test cases, the following dimensionless equations are now considered
\begin{subequations}
\begin{align}
\rho(\partial_t u + u \cdot \nabla u) + \nabla p &= \operatorname{N} \left( \nabla \times B \right) \times B + \frac{1}{\operatorname{Re}} \operatorname{div} \sigma + \frac{1}{\operatorname{Fr}} \rho \hat{g}, \\
\partial_t B -  \nabla \times \left(u \times B\right) &= - \frac{1}{\operatorname{Pm} \operatorname{Re}} \nabla \times \left( \nabla \times B \right), \\
\partial_t \rho + \nabla \cdot (\rho u) &= 0, \\
T (\partial_t s + \nabla \cdot (su)) & = \frac{1}{\operatorname{Re}} \sigma : \nabla u +  \frac{1}{\operatorname{Re} \operatorname{Pr}} \frac{\gamma}{\gamma-1} \Delta T + \frac{\operatorname{N}}{\operatorname{Pm} \operatorname{Re}} \left(\nabla \times B \right) \cdot \left(\nabla \times B \right),\\
\nabla \cdot B &= 0.
\end{align}
\end{subequations}
We are using standard notations, namely, $\operatorname{Re}$ is the Reynolds number, $\operatorname{N}$ the Stuart number, $\operatorname{Fr}$ the Froude number, $\operatorname{Pr}$ the Prandtl number, $\operatorname{Pm}$ the magnetic Prandtl number, $\gamma$ the heat capacity ratio, and $\hat{g}$ is a unit vector giving the orientation of the gravitational force. When the thermoelectric coupling is added, the coefficient $\alpha$ will be understood as a dimensionless number. Moreover, the density, pressure and temperature of reference are chosen such that the dimensionless perfect gas law simply reads
\begin{equation}
p = \rho T.
\end{equation}
Finally, the penalty parameter $\eta$ is set to $0.01$ in all computations.

\subsection{Reversible flow}

This first test case aims at evaluating the performances of the scheme in the absence of dissipation. No viscosity, resistivity or heat exchange are considered (i.e. $\operatorname{Re} = \operatorname{Pm} = \operatorname{Pr} = +\infty$). Gravitation is also removed ($\operatorname{Fr} = +\infty$) and $N$ is set to $0.014$. The heat capacity ratio is chosen as $\gamma = 1.4$.

\begin{figure}
\centering
\includegraphics[width = \textwidth]{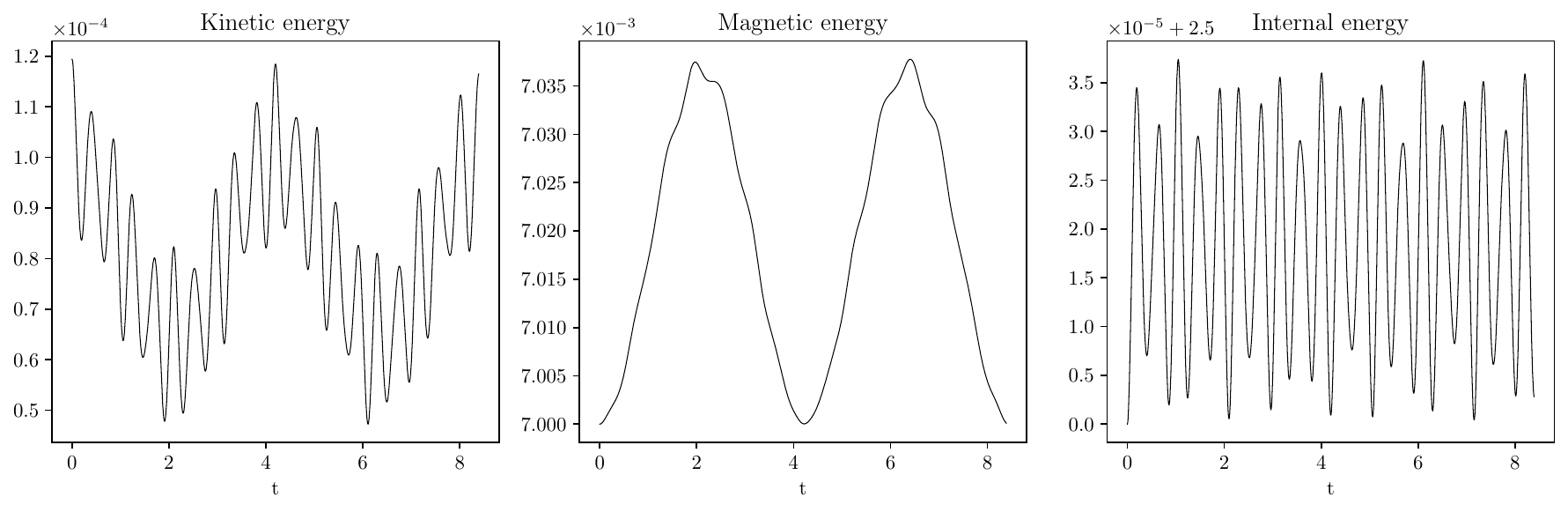}
\caption{Evolution of kinetic energy, internal energy and magnetic energy, all integrated over the whole domain. Both internal and magnetic energy evolves approximately as sinusoids and the kinetic energy is the superposition of their opposite. The ratio of frequencies is 10.}
\label{fig:reversible}
\end{figure}

The system under consideration is a square $[0,1] \times [0,1]$ with periodic boundary conditions. The density, temperature and magnetic field are initially homogeneous throughout the entire domain with $\rho = 1$, $T = 1$ and $B = (0, 1)$. A small velocity perturbation $(u_x,0)$, orthogonal to the magnetic field, is considered and reads
\begin{equation}
u_x = \left\{ 
\begin{array}{ll}
0.1 \exp\left( \frac{1}{(x - 0.5)^2 + (y - 0.5)^2 - 0.45^2} \right) & \text{if } (x - 0.5)^2 + (y - 0.5)^2 - 0.45^2 < 0, \\
0 & \text{otherwise}.
\end{array}
\right.
\end{equation}
This velocity compresses the fluid on the right, which in turn creates a pressure gradient opposing the fluid movement. Eventually, the fluid stops and is driven in the opposite direction by the pressure imbalance. This phenomenon is characterized by a transfer between kinetic energy and internal energy, not unlike the periodic transfer between kinetic energy and gravitational energy in the case of a swinging pendulum. Likewise, the initial velocity compresses magnetic field lines. The Lorentz force first slows the fluid down and then propels it in the opposite direction, thus leading to an alternate transfer between kinetic energy and magnetic energy. The system is then subject to two periodic phenomena, both characterized by their own frequency. For small deformations, their ratio should correspond to the ratio between the dimensionless speed of sound $\Tilde{v}_c$ and the dimensionless Alfv\'en speed $\Tilde{v}_A$
\begin{equation}
\frac{\Tilde{v}_c}{\Tilde{v}_A} = \sqrt{\frac{\gamma}{N}} = 10.
\end{equation}
Results for a $20\times20$ mesh, $r = s = 2$ and $\Delta t = 0.1$ are given on Figure \ref{fig:reversible}. The evolution of the different energies is consistent with the previous comment. In particular, the ratio of frequencies between the compressive and the magnetic modes is found to be 10. This test case showcases the good behavior of the scheme as oscillations are not visibly damped. This is critical as it means that numerical dissipation is basically absent and will not pollute physical dissipation when it is added.

\subsection{Two-dimensional magnetoconvection}

The scheme is now used to simulate Rayleigh--B\'enard convection in the presence of a magnetic field. In the case of an incompressible flow, the linear analysis of perturbations shows that stability is lost above a critical value of the Rayleigh number. It is here defined consistently with \cite{Hulburt1988} as
\begin{equation}
\operatorname{Ra} = \operatorname{Re}^2 (m+1) Z^2 \operatorname{Pr} (1 - (\gamma - 1)m)/\gamma.
\end{equation}
This threshold depends itself on the coupling between the fluid and the magnetic field, more specifically on the Chandrasekhar number
\begin{equation}
Q = \operatorname{N} \operatorname{Re}^2 \operatorname{Pm}.
\end{equation}
The coupling with a magnetic field tends to stabilize convection, so that the critical Rayleigh number increases with the value of $Q$. Although the present equations account for compressible effects and are highly non-linear, this general qualitative behavior will be shown to still be valid.

The domain $[0,2] \times [0,1]$ is considered with horizontal periodic boundary conditions. At the top and bottom, no-slip boundary conditions are set for the velocity and Dirichlet boundary conditions are enforced on the temperature $T_\text{bottom} = 1+Z$ and $T_\text{top} = 1$. Initial conditions follow that of \cite{Gawlik2024} with an additional vertical magnetic field $B = (0,1)$. The temperature is given a linear stratification $T = 1 + Z(1-z)$. Hydrostatic equilibrium then gives $\rho = T^m$ and $p = T^{m+1}$ with $\operatorname{Fr} = \frac{1}{(m+1) Z}$. An initial vertical velocity perturbation $u = (0, u_z)$ is considered with
\begin{equation}
u_z = \left\{ 
\begin{array}{ll}
\exp\left( \frac{1}{(x - 1)^2 + (y - 0.5)^2 - 0.2} \right) & \text{if } (x - 1)^2 + (y - 0.5)^2 - 0.2 < 0, \\
0 & \text{otherwise}.
\end{array}
\right.
\end{equation}
Dimensionless numbers are set to $\gamma = 0.1$, $m = 0$, $Z = 0.419524$ and $\operatorname{Pr} = \operatorname{Pm} = 2.5$. $\operatorname{Re}$ and $\operatorname{N}$ are not specified as they will take different values depending on the Rayleigh and Chandrasekhar numbers. Finally, simulations are performed on a $32 \times 16$ mesh with $r = s = 1$ and a time step of $0.1$.

\begin{figure}
\centering
\includegraphics[width = \textwidth]{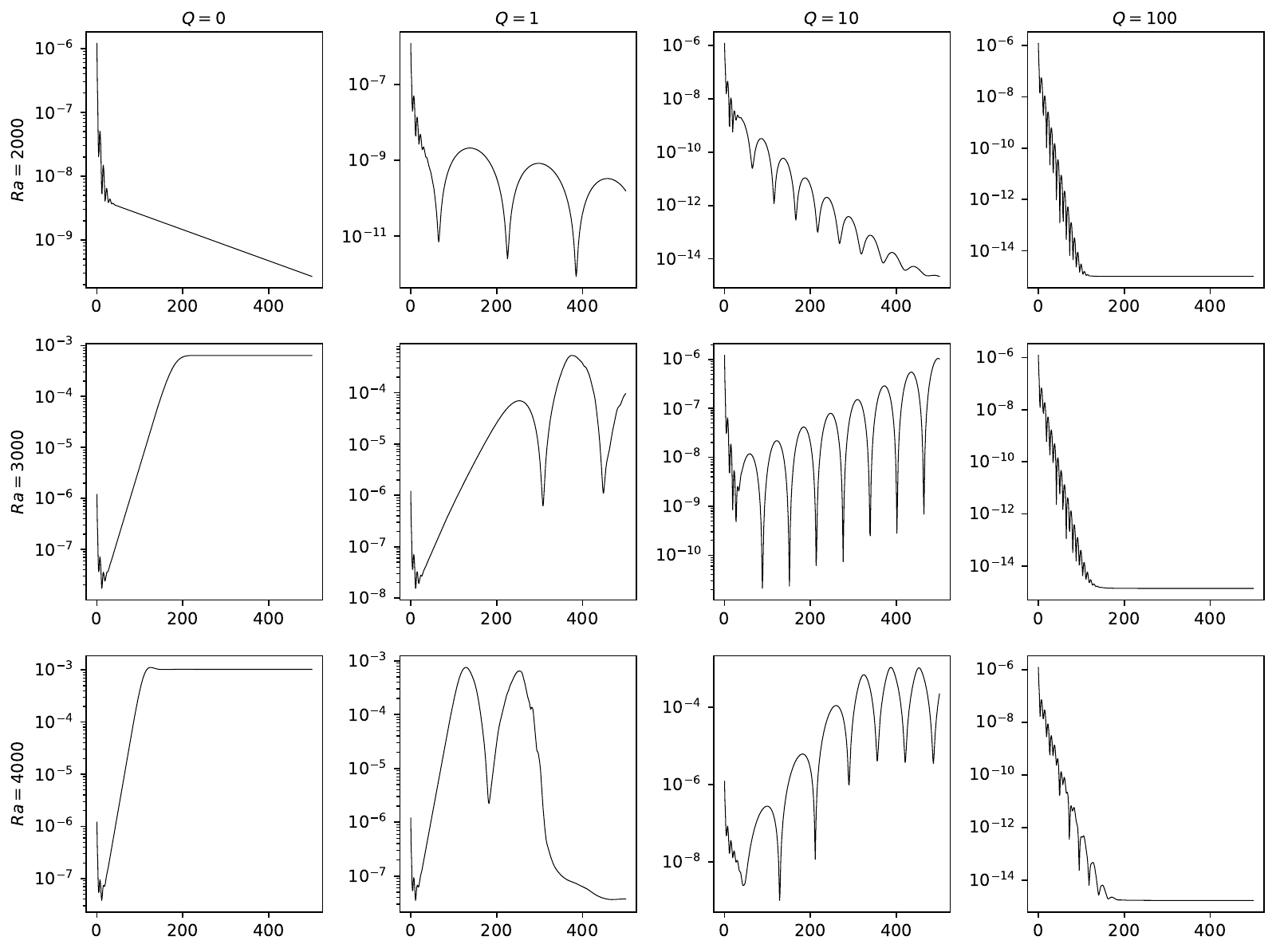}
\caption{Evolution of the kinetic energy for different values of $\operatorname{Ra}$ and $\operatorname{Q}$. Instability is associated with strong $\operatorname{Ra}$ and weak $\operatorname{Q}$.}
\label{fig:RB_stability}
\end{figure}

\begin{figure}
\centering
\includegraphics[width = \textwidth]{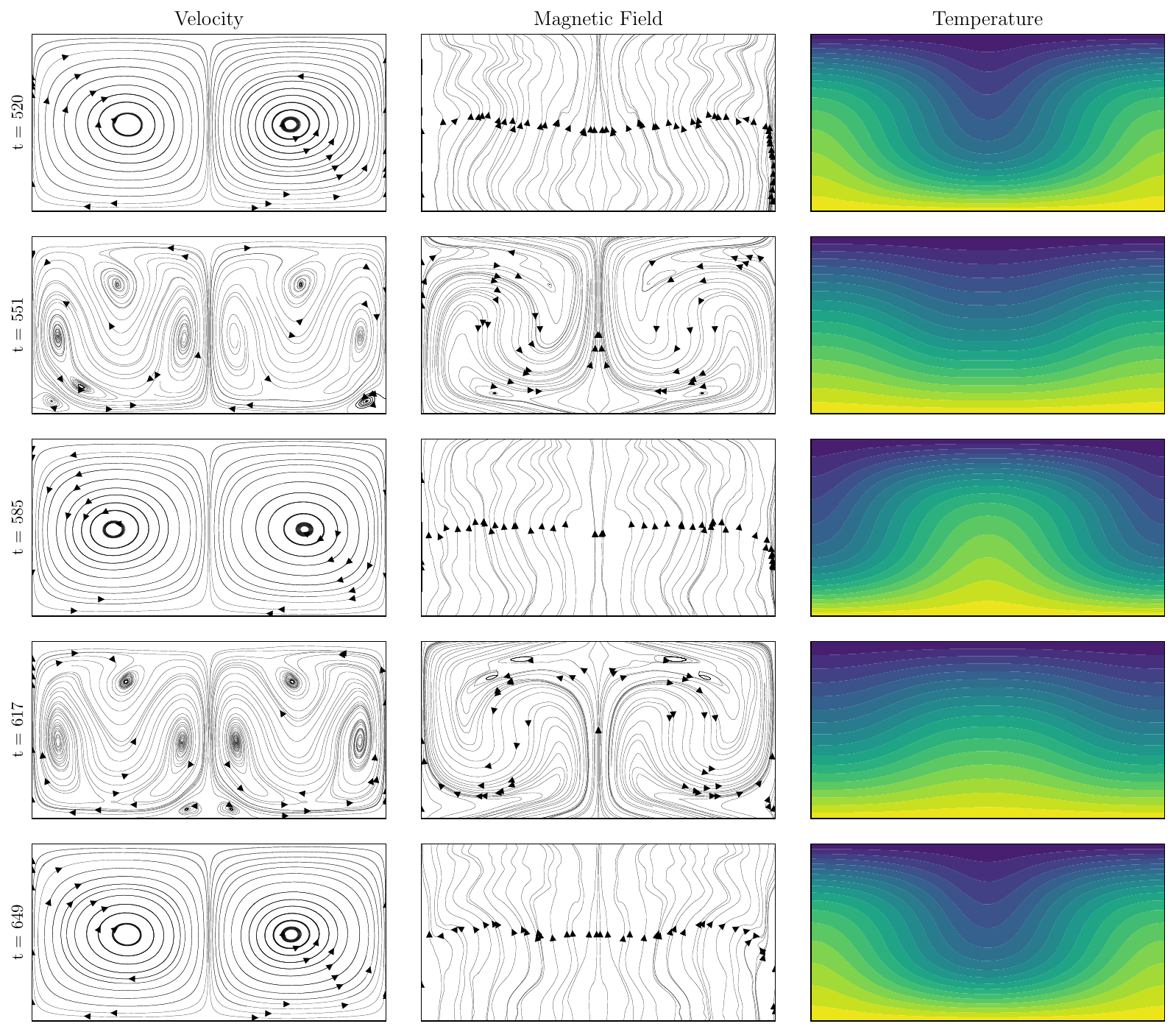}
\caption{Snapshots of a periodic unstable flow ($\operatorname{Ra} = 4000$ and $\operatorname{Q} = 10$) at different times. It shows the reversal of the convection rolls due to the deformation of the magnetic field lines.}
\label{fig:RB_periodic}
\end{figure}

\paragraph{Flow stability.} The evolution of the system kinetic energy for different values of $\operatorname{Ra}$ and $\operatorname{Q}$ is shown in Figure \ref{fig:RB_stability}. As expected, it is seen that increasing $Q$ slows down the onset of convection or suppresses it altogether. With $Q = 100$, kinetic energy quickly reaches machine precision with the three values of $\operatorname{Ra}$.
Contrary to classical Rayleigh--B\'enard instability (i.e. without the contribution of a magnetic field), convection is not stationary. Rather, it displays a periodic behavior translating to oscillations on the kinetic energy (see Figure \ref{fig:RB_stability}). Phenomenologically, these oscillations are similar to those of the magnetic energy in the previous test case: As the magnetic field is distorted by the fluid movement, kinetic energy is transferred to magnetic energy. Once enough magnetic energy has been accumulated, it is converted back to kinetic energy through the inversion of the convection rolls. Hence, the magnetic field has the effect of an axial spring on the fluid movement. This phenomenon is illustrated in Figure \ref{fig:RB_periodic}.

\begin{figure}
\centering
\includegraphics[width = \textwidth]{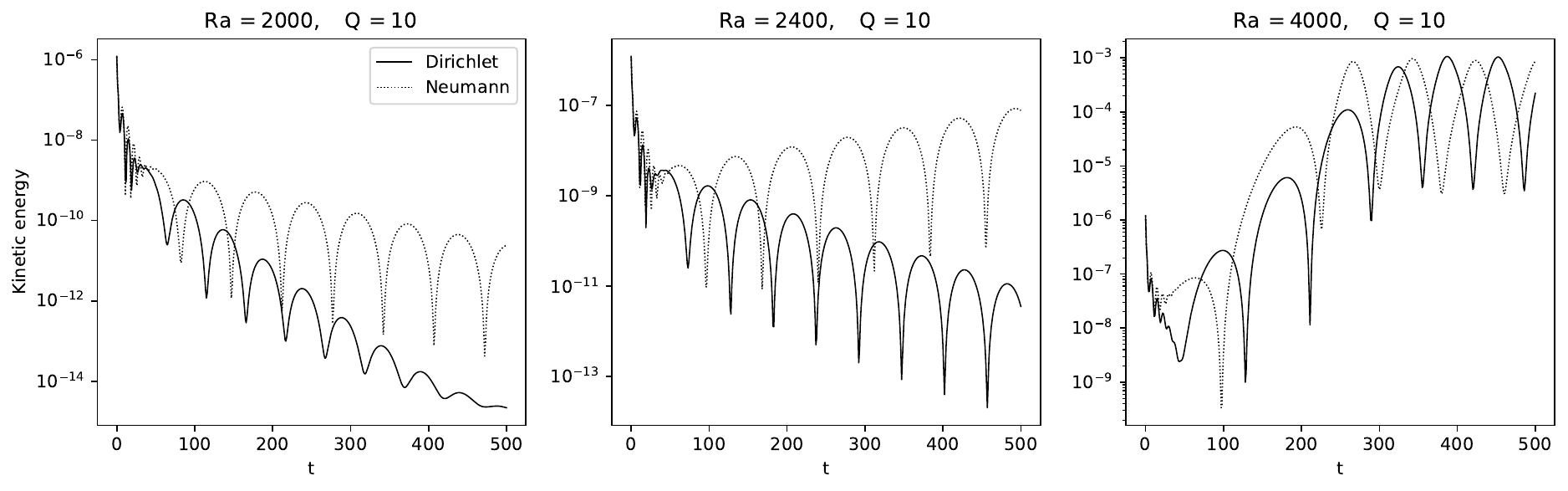}
\caption{Evolution of the kinetic energy over time for different values of $\operatorname{Ra}$ and different boundary conditions. Neumann boundary conditions are more conducive to the onset of convection than Dirichlet boundary conditions. Reversal of the flow also seems to occur at smaller frequencies.}
\label{fig:RB_neumann}
\end{figure}

\paragraph{Impact of boundary conditions.} Results with Neumann boundary conditions are shown in Figure \ref{fig:RB_neumann}. They seem to induce more instability than Dirichlet ones. For $\operatorname{Ra} = 2000$, both give a stable flow, but the kinetic energy decreases faster with Dirichlet boundary conditions. For $\operatorname{Ra} = 2400$, Dirichlet boundary conditions induce a stable flow while Neumann ones do not. Changes in boundary conditions also slightly affect the frequency of flow reversals.

\begin{figure}
\centering
\includegraphics[width = \textwidth]{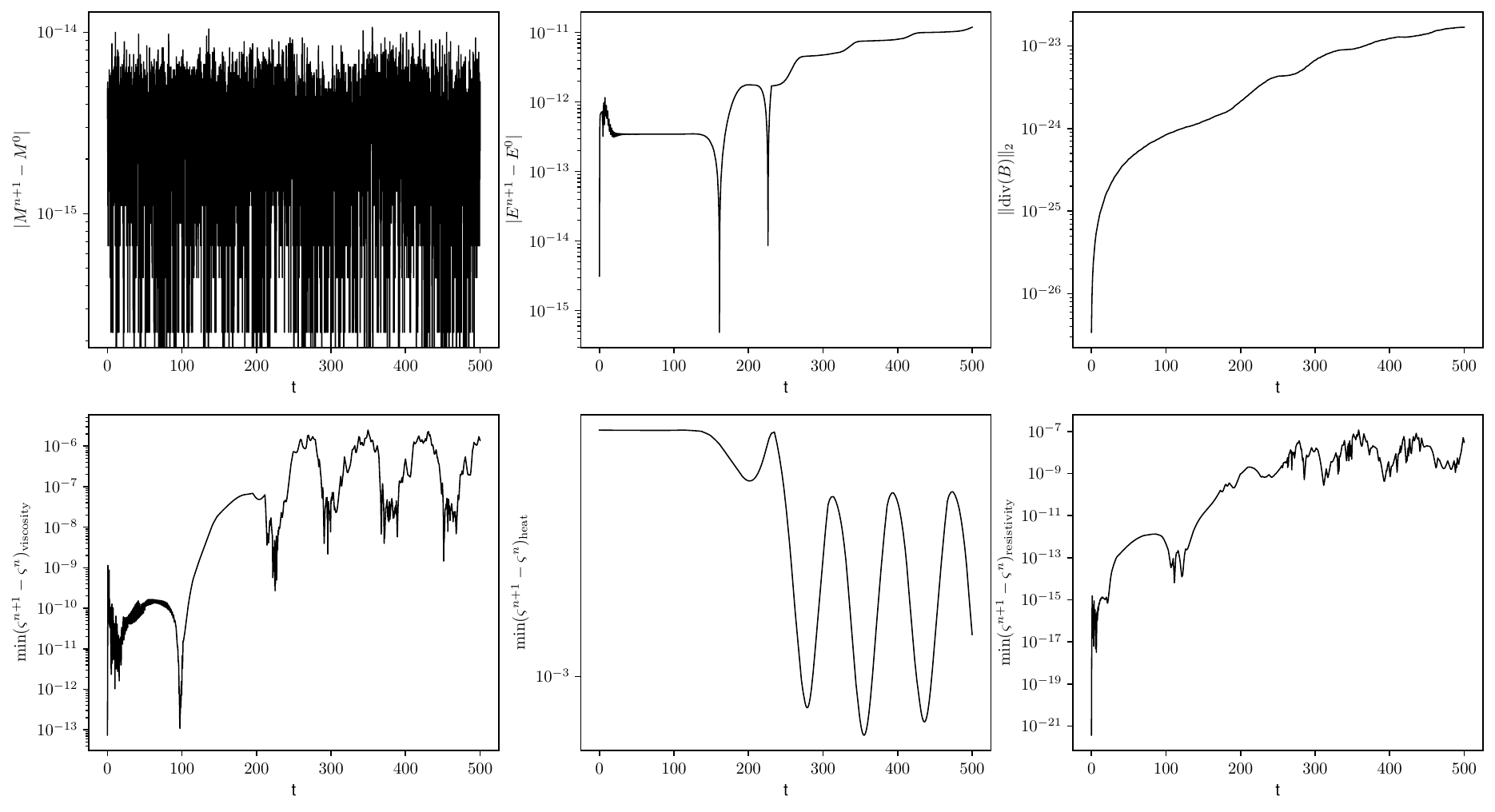}
\caption{The properties of the scheme are assessed on the convection test case with Neumann boundary convection and $(\operatorname{Ra}, \operatorname{Q}) = (4000, 10)$. On the second row is plotted the minimum value of the contribution of the different sources of entropy production. Compliance with mass conservation, Gauss' magnetic law and the two principles of thermodynamics is confirmed at the discrete level.}
\label{fig:conservation_neumann}
\end{figure}

\paragraph{Conservation and balance laws.} With the previous Neumann boundary conditions, the top and bottom energy fluxes balance one another so that the total energy is conserved (even though the system is not isolated). The properties of the scheme are then assessed for this test case and results can be found in Figure \ref{fig:conservation_neumann}. It is shown that mass and energy are both conserved at the discrete level. So is Gauss' magnetic law. Finally, the different sources of entropy production are shown to individually have a positive contribution to the total entropy variation.

\paragraph{Strong deformations of the magnetic field.} Finally, it is important to emphasize that the more inertia there is (high $\operatorname{Ra}$) or the weaker the coupling is (low $\operatorname{Q}$), the more the magnetic field lines are distorted by the flow. This can lead to two issues. First, strong deformations can lead to oscillations, thus jeopardizing the scheme's stability. Second, if the dimension of the finite element space is not big enough to capture these deformations, high frequency magnetic energy will instead populate low frequency levels and create a magnetic smudge, thus leading to some questionable numerical results. This is what happens in Figure \ref{fig:RB_stability} for $(\operatorname{Ra}, \operatorname{Q}) = (4000, 1)$ for which periodic modes cannot properly set. While the first issue could be solved with some form of artificial resistivity \cite{Fambri2025}, the second is a mathematical limit that can only be overcome by adding a sufficient number of elements.

\subsection{Magnetoconvection and thermoelectric effect}

No thermoelectric effect can exist in two dimensions when the magnetic field is parallel to the plane. It is however possible to include the coupling when the magnetic field is orthogonal. The present test case takes the initial conditions of the previous one, with Dirichlet boundary conditions on the temperature, while considering a magnetic field $B = (0, 0, 1)$.

\begin{figure}
\centering
\includegraphics[width = 0.4\textwidth]{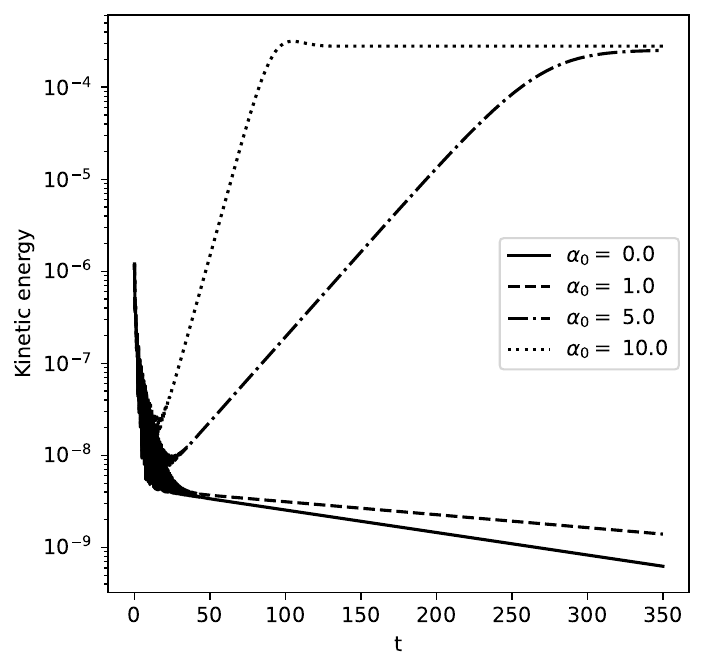}
\caption{Evolution of the kinetic energy for different values of coefficient $\alpha_0$.}
\label{fig:orthogonal_convection}
\end{figure}

$\operatorname{Ra}$ and $\operatorname{Q}$ are respectively taken as $2000$ and $10$. Finally, simulations are performed on a $32 \times 16$ mesh with $r = s = 1$ and a time step of $0.1$. The (dimensionless) thermoelectric coefficient is here given as a function of the density and temperature
\begin{equation}
\alpha(\rho, T) = \alpha_0 \frac{T}{\rho}.
\end{equation}
Results with and without thermoelectric effect are shown on Figure \ref{fig:orthogonal_convection}. While the flow is stable without the thermoelectric effect, it becomes unstable for a sufficiently large value of $\alpha_0$. Besides, the larger $\alpha_0$ is, the more kinetic energy the system has when it reaches a stationary flow, and the less time it takes to reach that state.

\section{Conclusion}

A numerical method for visco-resistive MHD with heat transfer and thermoelectric effect was developed based on the weak form of the equations, derived from a variational formulation in nonequilibrium thermodynamics. Hamilton's principle, extended to systems with dissipation, is shown to be a strong source of guidance to derive schemes which preserve important features of these complex flows. In particular, the two laws of thermodynamics are satisfied at the fully discrete level, with energy being preserved and entropy not allowed to be destroyed, for various thermal boundary conditions. Numerical results confirms the accuracy and stability of the scheme while showcasing how dissipation modifies MHD flows.

The present work paves the way for further studies. The same approach could be used for the numerical approximation of other equations that can be derived from a variational principle. An important example is that of multiphase flows which have a large spectrum of applications in theoretical science and industry alike. Such equations also allow for a virtually infinite number of dissipative processes and couplings, thus making consistency with thermodynamics a critical and challenging target. On the other hand, the present method could be adapted to deal with challenging regimes including stiff gradients or general equations of states. \textit{Ad hoc} numerical modifications would then be needed to ensure positive temperature and densities, as well as accurate and robust shock-capturing.

\paragraph{Acknowledgment.} We are thankful for fruitful discussions with Xavier Garbet. Evan Gawlik was supported by NSF grant DMS-2533499 and the Simons Foundation award MPS-TSM-00002615. Fran\c{c}ois Gay-Balmaz was  supported by a startup grant from Nanyang Technological University. Fran\c{c}ois Gay-Balmaz and Bastien Manach-P\'erennou were supported by the National
Research Foundation Singapore (NRF) core funding ‘Fusion Science for Clean Energy’.

\bibliographystyle{plain}
\bibliography{bib.bib}

\appendix

\section{Consistency of the discrete Seebeck-Peltier coupling}
\label{app:consistency}

To show that this choice of $h_h$ gives rise to a consistent discretization of the equations of motion, let $H$ be a smooth vector field, let $T$ and $\alpha$ be smooth scalar fields, and let $w$ be piecewise smooth but discontinuous. Then
\begin{subequations}
\begin{align}
h_h(1,Tw ,H) 
& = -\int_{\partial\Omega} H \cdot (\alpha n \times \nabla (Tw)) \, {\rm d} a +  \sum_{K \in \mathcal{T}_h} \int_K H \cdot (\nabla \alpha \times \nabla (Tw)) \, {\rm d} x \nonumber \\
\label{1HTw}
& \quad - \sum_{e \in \mathcal{E}_h^0} \int_e H \cdot ( \nabla \alpha \times T \llbracket w \rrbracket) \, {\rm d} a, \\
h_h(w,T,H)
& = -\int_{\partial\Omega} H \cdot (w\alpha n \times \nabla T) \, {\rm d} a +  \sum_{K \in \mathcal{T}_h} \int_K H \cdot (\nabla (w\alpha) \times \nabla T) \, {\rm d} x \nonumber \\
\label{wHT}
&\quad - \sum_{e \in \mathcal{E}_h^0} \int_e H \cdot (\llbracket w \rrbracket \alpha \times \nabla T ) \, {\rm d} a.
\end{align}
\end{subequations}
Subtracting~\eqref{wHT} from~\eqref{1HTw}, and using the identities
\begin{subequations}
\begin{align}
-\alpha n \times \nabla (Tw) + w\alpha n \times \nabla T & = -T\alpha n \times \nabla w, \\
\nabla \alpha \times \nabla (Tw) - \nabla (w\alpha) \times \nabla T & = \nabla (T\alpha) \times \nabla w, \\
-\nabla \alpha \times T\llbracket w \rrbracket + \llbracket w \rrbracket \alpha \times \nabla T & = -\nabla(T\alpha) \times \llbracket w \rrbracket,
\end{align}
\end{subequations}
it yields
\begin{multline}
h_h(1,Tw,H) - h_h(w,T,H) = -\int_{\partial\Omega} H \cdot (T\alpha n \times \nabla w) \, {\rm d} a  \\ + \sum_{K \in \mathcal{T}_h} \int_K (H \times \nabla (T\alpha)) \cdot \nabla w \, {\rm d} x - \sum_{e \in \mathcal{E}_h^0} \int_e (H \times \nabla(T\alpha)) \cdot \llbracket w \rrbracket \, {\rm d} a.
\end{multline}
We can rewrite the integrals over $e \in \mathcal{E}_h^0$ as integrals over boundaries of $n$-simplices, bearing in mind that edges in $\partial\Omega$ are not included in the set $\mathcal{E}_h^0$.  We obtain
\begin{align*}
h_h(1,Tw,H) - h_h(w,T,H)
&= -\int_{\partial\Omega} H \cdot (T\alpha n \times \nabla w) \, {\rm d} a + \sum_{K \in \mathcal{T}_h} \int_K (H \times \nabla (T\alpha)) \cdot \nabla w \, {\rm d} x \nonumber \\
&\quad\quad - \sum_{K \in \mathcal{T}_h} \int_{\partial K} (H \times \nabla(T\alpha)) \cdot wn \, {\rm d} a + \int_{\partial \Omega} (H \times \nabla(T\alpha)) \cdot wn \, {\rm d} a.
\end{align*}
The integrals over $K \in \mathcal{T}_h$ can be combined using integration by parts, and the integrals over $\partial\Omega$ can be combined as well, resulting in
\begin{align*}
h_h(1,Tw,H) - h_h(w,T,H)
&= -\int_\Omega \operatorname{div}(H \times \nabla(T\alpha)) w \, {\rm d} x + \int_{\partial\Omega} (H \times \nabla (wT\alpha)) \cdot n \, {\rm d} a \\
&= -\int_\Omega \curl H \cdot \nabla(T\alpha) w \, {\rm d} x + \int_{\partial\Omega} (H \times \nabla (wT\alpha)) \cdot n \, {\rm d} a.
\end{align*}
The second term above equals $\int_{\partial\Omega} wT\alpha \curl H \cdot n \, {\rm d} a$ because
\begin{subequations}
\begin{align}
0 &= \int_\Omega \dv(\curl(wT\alpha H)) \, {\rm d} x \\
&= \int_{\partial\Omega} \curl(wT\alpha H) \cdot n \, {\rm d} a \\
&= \int_{\partial\Omega} (\nabla (wT\alpha) \times H) \cdot n \, {\rm d} a + \int_{\partial\Omega} wT\alpha \curl H \cdot n \, {\rm d} a.
\end{align}
\end{subequations}
Thus,
\begin{equation}
\label{eq:consistency_h}
h_h(1,Tw,H) - h_h(w,T,H) = -\int_\Omega \curl H \cdot \nabla(T\alpha) w \, {\rm d} x + \int_{\partial\Omega} wT\alpha \curl H \cdot n \, {\rm d} a.
\end{equation}
Furthermore, it is shown in \cite[§3.3]{Gawlik2024} that, for smooth $T$, $d_h$ satisfies
\begin{equation}
\label{eq:consistency_d}
d_h(1,T,Tw) - d_h(w,T,T)
= \int_\Omega \kappa \Delta T w \, {\rm d} x - \int_{\partial\Omega} \kappa \nabla T  \cdot n w \, {\rm d} a.
\end{equation}
Combining \eqref{eq:consistency_h} and \eqref{eq:consistency_d} eventually yields
\begin{equation}
\overline{d}_h(1,T,Tw,H) - \overline{d}_h(w,T,T,H)
= -\int_\Omega \dv(Tj_s) w \, {\rm d} x + \int_{\partial\Omega} Tj_s \cdot n w \, {\rm d} a.
\end{equation}
The discretization is thus consistent and the boundary condition $j_s \cdot n=0$ is enforced weakly.

\end{document}